\documentclass[11pt]{amsart}
\usepackage{amsmath,amssymb,amsthm,mathtools,enumitem,booktabs}
\usepackage[margin=2.7cm]{geometry}
\usepackage{tikz}
\usetikzlibrary{arrows.meta,calc}
\usepackage[colorlinks=true,linkcolor=blue,citecolor=blue,urlcolor=blue]{hyperref}

\newtheorem{theorem}{Theorem}[section]
\newtheorem{lemma}[theorem]{Lemma}
\newtheorem{proposition}[theorem]{Proposition}
\newtheorem{corollary}[theorem]{Corollary}
\newtheorem{conjecture}[theorem]{Conjecture}
\theoremstyle{definition}
\newtheorem{definition}[theorem]{Definition}
\newtheorem{example}[theorem]{Example}
\newtheorem{question}[theorem]{Question}
\theoremstyle{remark}
\newtheorem{remark}[theorem]{Remark}

\numberwithin{equation}{section}

\newcommand{\AT}{\operatorname{AT}}
\newcommand{\ch}{\operatorname{ch}}
\newcommand{\CT}{\operatorname{CT}}
\newcommand{\sgn}{\operatorname{sgn}}
\newcommand{\coef}[1]{\bigl[#1\bigr]}
\newcommand{\one}{\mathbf 1}
\newcommand{\Q}{\mathbb{Q}}
\newcommand{\N}{\mathbb{N}}

\title{An Alon--Tarsi analogue of Ohba's theorem}
\author{Witold Sikora}
\address{Faculty of Mathematics and Computer Science, Jagiellonian University, Krak\'ow, Poland}
\email{witold.sikora@student.uj.edu.pl}

\begin{document}

\begin{abstract}\emergencystretch=2em
Ohba proved that joining a graph $G$ to a large enough clique $K_k$ produces a chromatic-choosable graph. We study the same question for the Alon--Tarsi number: let $f(G)$ be the least $k$ with $\AT(G\vee K_k)=\chi(G)+k$. We show that $f(G)$ is always finite, that $f(G)\le n^2/8$ for every $n$-vertex graph, tight up to an additive $O(n)$, and that with high probability $f(G(n,1/2))=(\frac12+o(1))\,n\log_2n$. Thus quadratically many universal vertices may be needed, as for DP-colouring and unlike list colouring. We also show that every graph with at most $\chi(G)+2\sqrt{2\chi(G)-2}-2$ vertices has Alon--Tarsi number $\chi(G)$; this is sharp up to an additive constant and improves a bound of Kozik, Micek and Zhu for the paint number.

Our main structural result is that, for a uniquely $\chi$-colourable graph $G$, the whole family $G\vee K_k$ is governed by a single multigraph $M_G$ on $\chi$ vertices: $f(G)=\AT(M_G)-\chi$. Combined with the constant term identities of Dyson and Morris, this gives exact values of $f$ for bipartite graphs and for complete multipartite graphs in which all parts but one have the same size. For complete tripartite graphs we determine $f$ up to an additive error of at most one. The proofs remove the universal vertices with the dual Cauchy identity for Schur functions and then use Lagrange interpolation on the colour classes.
\end{abstract}

\maketitle

\section{Introduction}\label{sec:intro}

Throughout, graphs are finite and simple and have at least one vertex. The \emph{join} $G\vee H$ of two vertex-disjoint graphs is obtained from their disjoint union by adding all edges between $V(G)$ and $V(H)$. Thus $G\vee K_k$ is $G$ together with $k$ new vertices, called \emph{universal vertices}, which are adjacent to each other and to every vertex of $G$.

For a graph $G$ with vertices $v_1<\dots<v_n$, its \emph{graph polynomial} is
\[
P_G(x)=\prod_{\substack{v_iv_j\in E(G)\\ i<j}}(x_{v_i}-x_{v_j}).
\]
The \emph{Alon--Tarsi number} $\AT(G)$ is the least $d$ such that $P_G$ contains a monomial with nonzero coefficient in which every exponent is at most $d-1$. For a loopless multigraph $M$, define $P_M$ by the same product, with each parallel edge contributing a repeated factor; define $\AT(M)$ from $P_M$ in the same way. The point of this definition is Alon's Combinatorial Nullstellensatz \cite{AlonCN}: if the coefficient of $\prod_vx_v^{a_v}$ in $P_G$ is nonzero, then $G$ can be properly coloured from any lists of sizes $a_v+1$. Hence $\AT(G)$ bounds the \emph{choice number} $\ch(G)$, the least $d$ such that $G$ can be properly coloured from arbitrary lists of $d$ colours \cite{AT}. Schauz \cite{SchauzFlex} showed that it also bounds the \emph{paint number} $\chi_P(G)$, the online version of $\ch(G)$ introduced in \cite{SchauzPaint}, in which the lists are revealed one colour at a time. Altogether,
\[
\chi(G)\le\ch(G)\le\chi_P(G)\le\AT(G).
\]

A graph is \emph{chromatic-choosable} if $\ch(G)=\chi(G)$. Ohba \cite{Ohba} proved that $G\vee K_k$ is chromatic-choosable for every graph $G$ and all sufficiently large $k$; by the theorem of Noel, Reed and Wu \cite{NRW}, $k\ge|V(G)|-2\chi(G)-1$ suffices. The least such $k$ is the \emph{Ohba number} of $G$. The same question has been studied for the paint number \cite{HWZ,KKLZ,KMZ} and for the DP-chromatic number $\chi_{DP}$, a generalisation of the choice number in which, along each edge, the colours of the two endpoints that conflict are prescribed by an arbitrary matching rather than by equality. For $\chi_{DP}$ the answer is different. Let $Z_{DP}(G)$ be the least $k$ with $\chi_{DP}(G\vee K_k)=\chi(G)+k$, and $Z_{DP}(n)$ its maximum over $n$-vertex graphs. Bernshteyn, Kostochka and Zhu \cite{BKZ} showed that quadratically many universal vertices may be needed, and Zhang and Dong \cite{ZD} that $\lceil2.4\,e(G)\rceil$ always suffice, where $e(G)$ is the number of edges. Altogether $n^2/4-O(n)\le Z_{DP}(n)\le n^2-(n+3)/2$.

This paper is about the Alon--Tarsi number. Since $\chi(G\vee K_k)=\chi(G)+k$, the natural quantity is
\[
f(G)=\min\{k\ge0:\AT(G\vee K_k)=\chi(G)+k\}.
\]
The gap $\AT(G\vee K_k)-\chi(G)-k$ does not increase with $k$ (Lemma~\ref{lem:gap}), but there is no obvious reason for it to vanish: the colouring number (one more than the degeneracy) of $G\vee K_k$, for instance, exceeds $\chi(G)+k$ by an amount independent of $k$. Nor do the results for choosability or the paint number imply that $f(G)$ is finite, since $\AT$ dominates both parameters. We write $\chi=\chi(G)$ and $n=|V(G)|$; further notation is collected in Section~\ref{sec:prelim}.

\subsection{How many universal vertices are needed?}
Our first result shows that $f(G)$ is finite and determines its order of magnitude.

\begin{theorem}\label{thm:main}
Let $G$ be a graph with a proper $\chi$-colouring in which every colour class has at most $m$ vertices. Then
\[
\max_{G'\subseteq G}\frac{e(G')-(\chi-1)|V(G')|}{\chi}\ \le\ f(G)\ \le\ (\chi-1)\Bigl\lfloor\frac{(m-1)^2}{2}\Bigr\rfloor.
\]
The upper bound is attained by the complete $\chi$-partite graph $K_{m\star\chi}$ with parts of size $m$, and for even $m$ so is the lower bound.
\end{theorem}

The upper bound is efficient for balanced colourings. A second bound depends only on the number of edges leaving each colour class, and it settles the extremal function $F(n)=\max\{f(G):|V(G)|=n\}$ asymptotically.

\begin{theorem}\label{thm:intro-edges}
For every proper colouring of $G$ with $\chi=\chi(G)$ classes $V_1,\dots,V_\chi$,
\[
f(G)\le\max_i\Bigl\lfloor\frac{e(V_i,V(G)\setminus V_i)}{2}\Bigr\rfloor .
\]
Together with Theorem~\ref{thm:bipartite}, this gives $n^2/8-n/2\le F(n)\le n^2/8$. Moreover, for the uniform random graph $G(n,1/2)$, with probability tending to $1$ as $n\to\infty$,
\[
f(G(n,1/2))=\bigl(\tfrac12+o(1)\bigr)\,n\log_2n .
\]
\end{theorem}

Since $Z_{DP}(n)\ge n^2/4-O(n)$, in the worst case the Alon--Tarsi number needs at most about half as many universal vertices as the DP-chromatic number. The balanced complete bipartite graph is asymptotically extremal. In fact, $f(K_{m,m})=\lfloor(m-1)^2/2\rfloor$, whereas $\ch(K_{m,m}\vee K_k)=k+2$ already for $k\ge2m-5$ \cite{NRW}; thus the Alon--Tarsi number of $K_{m,m}\vee K_k$ exceeds its choice number throughout a range of $k$ whose length is quadratic in $m$.

The same bounds answer a question about graphs with few vertices. Every graph $G$ with $|V(G)|\le\chi(G)+2\sqrt{2\chi(G)-2}-2$ satisfies $\AT(G)=\chi(G)$ (Theorem~\ref{thm:few}), and the graphs $K_{m,m}\vee K_k$ show that this is best possible up to an additive constant: the least number of vertices of a $\chi$-chromatic graph with $\AT(G)>\chi$ is $\chi+2\sqrt{2\chi}+O(1)$ (Corollary~\ref{cor:psi}). Since $\chi_P\le\AT$, this improves the bound $\chi+\sqrt\chi$ of Kozik, Micek and Zhu \cite{KMZ} towards the on-line version of Ohba's conjecture, which predicts $\chi_P(G)=\chi(G)$ whenever $|V(G)|\le2\chi(G)$. It also shows that the Alon--Tarsi number alone cannot prove that conjecture.

\subsection{Collapse to a multigraph}
The main result is an exact reduction for uniquely colourable graphs. Here \emph{uniquely $\chi$-colourable} means that $V(G)$ has exactly one partition into $\chi$ independent sets. Let $G$ be uniquely $\chi$-colourable with colour classes $V_1,\dots,V_\chi$ of sizes $m_1,\dots,m_\chi$, and let $M_G$ be the multigraph on $\{1,\dots,\chi\}$ in which $i$ and $j$ are joined by
\begin{equation}\label{eq:mu}
\mu_{ij}=e(V_i,V_j)-m_i-m_j+2
\end{equation}
parallel edges. Here $e(V_i,V_j)$ is the number of edges between $V_i$ and $V_j$. Thus $\mu_{ij}-1$ is the circuit rank (the number of independent cycles) of the bipartite graph $G[V_i\cup V_j]$ induced by $V_i\cup V_j$, which is connected because $G$ is uniquely colourable (if it were disconnected, exchanging the colours $i$ and $j$ on one component would give a second colouring); in particular $\mu_{ij}\ge1$. For the complete multipartite graph $K_{m_1,\dots,m_\chi}$ one has $\mu_{ij}=(m_i-1)(m_j-1)+1$.

\begin{theorem}\label{thm:collapse}
If $G$ is uniquely $\chi$-colourable, then for every $k\ge0$
\[
\AT(G\vee K_k)=\chi+k\quad\Longleftrightarrow\quad\AT(M_G)\le\chi+k.
\]
Consequently $f(G)=\AT(M_G)-\chi$.
\end{theorem}

Here $\AT(M_G)\ge\chi$, because $\mu_{ij}\ge1$ for all $i<j$, so that $M_G$ contains $K_\chi$.

\begin{figure}[ht]
\centering
\begin{tikzpicture}[scale=0.9,
  v/.style={circle,draw,inner sep=1.6pt,fill=#1},
  u/.style={rectangle,draw,inner sep=2.2pt,fill=black!15}]
\begin{scope}
\node[v=red!55]  (a1) at (-0.45,2.2) {}; \node[v=red!55] (a2) at (0.45,2.2) {};
\node[v=blue!45] (b1) at (-1.6,0) {}; \node[v=blue!45] (b2) at (-1.0,-0.6) {};
\node[v=green!50!black] (c1) at (1.0,-0.6) {}; \node[v=green!50!black] (c2) at (1.6,0) {}; \node[v=green!50!black] (c3) at (2.0,0.7) {};
\foreach \x in {a1,a2}{\foreach \y in {b1,b2,c1,c2,c3}{\draw[black!40] (\x)--(\y);}}
\foreach \x in {b1,b2}{\foreach \y in {c1,c2,c3}{\draw[black!40] (\x)--(\y);}}
\node[red!70!black] at (0,2.65) {\scriptsize $V_1$};
\node[blue!70!black] at (-1.85,-0.55) {\scriptsize $V_2$};
\node[green!40!black] at (2.15,-0.45) {\scriptsize $V_3$};
\draw[dashed,black!50] (0.2,0.75) ellipse (2.6 and 2.25);
\node at (0.2,-1.85) {$G=K_{2,2,3}$};
\end{scope}
\begin{scope}[shift={(4.9,0.75)}]
\node[u] (y1) at (0,0.9) {}; \node[u] (y2) at (0,0.3) {}; \node at (0,-0.15) {\scriptsize$\vdots$}; \node[u] (y3) at (0,-0.75) {};
\draw[rounded corners=4pt,black!60] (-0.45,-1.2) rectangle (0.45,1.35);
\node at (0,-2.6) {$K_k$};
\end{scope}
\draw[line width=2.2pt,black!25] (2.85,0.75) -- (4.4,0.75);
\node[align=center] at (3.6,1.25) {\scriptsize join};
\draw[-{Stealth},thick] (5.8,0.75) -- node[above] {\small collapse} (7.4,0.75);
\begin{scope}[shift={(8.3,-0.2)}]
\node[v=red!55,label=above:{\small $1$}] (A) at (0.9,2.3) {};
\node[v=blue!45,label=left:{\small $2$}] (B) at (0,0) {};
\node[v=green!50!black,label=right:{\small $3$}] (C) at (1.8,0) {};
\draw (A) to[bend left=10] (B); \draw (A) to[bend right=10] (B);
\draw (A) to[bend left=16] (C); \draw (A) -- (C); \draw (A) to[bend right=16] (C);
\draw (B) to[bend left=16] (C); \draw (B) -- (C); \draw (B) to[bend right=16] (C);
\node at (0.05,1.35) {\scriptsize $2$}; \node at (1.75,1.35) {\scriptsize $3$}; \node at (0.9,-0.6) {\scriptsize $3$};
\node at (0.9,-1.65) {$M_G$};
\end{scope}
\end{tikzpicture}
\caption{Theorem~\ref{thm:collapse} for $G=K_{2,2,3}$. The $k$ universal vertices (squares) are joined to everything. The multigraph $M_G$ has multiplicities $\mu_{ij}=(m_i-1)(m_j-1)+1$. Here $P_{M_G}=(u_1-u_2)^2(u_1-u_3)^3(u_2-u_3)^3$ has degree $8$, so every monomial has an exponent at least $3$, and the coefficient of $u_1^3u_2^2u_3^3$ is $8\ne0$. Hence $\AT(M_G)=4$, so $\AT(G\vee K_k)=3+k$ exactly when $k\ge1$, and $f(G)=1$.}\label{fig:collapse}
\end{figure}
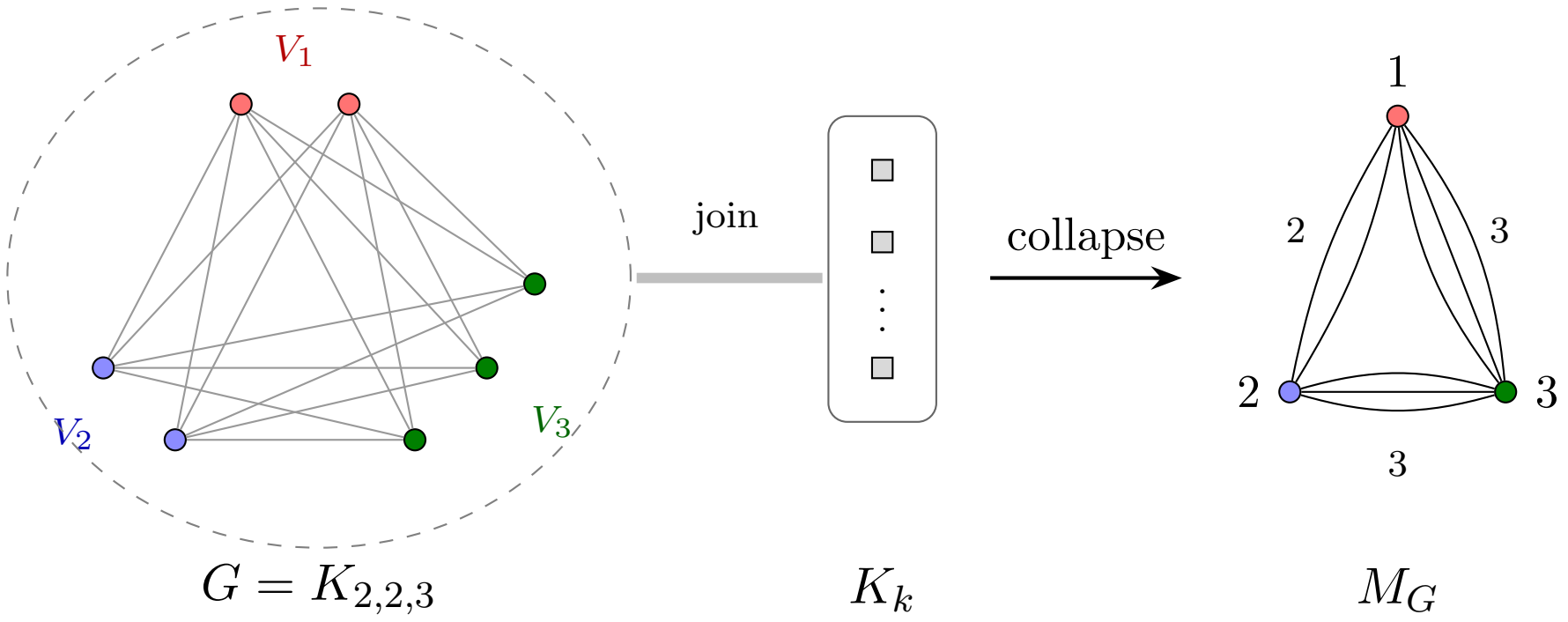

Theorem~\ref{thm:collapse} is the central result of the paper. A priori, each $k$ poses a separate problem: the graph polynomial of $G\vee K_k$ has $n+k$ variables and degree $e(G)+nk+\binom k2$, and the monomials that have to be examined change with $k$. Theorem~\ref{thm:collapse} replaces this infinite family of problems by a single one, about a multigraph $M_G$ that does not depend on $k$ and has $\chi$ vertices rather than $n$ (Figure~\ref{fig:collapse}). The universal vertices disappear completely, and the structure of $G$ enters only through the $\binom\chi2$ numbers $\mu_{ij}$. In particular $f(G)$ is determined by one Alon--Tarsi number of a small multigraph. This is what makes exact values accessible: for the complete multipartite graphs treated below, $f$ is found by evaluating or bounding coefficients of the single polynomial $P_{M_G}$, by a constant term identity (Theorems~\ref{thm:equal} and \ref{thm:morris}), by a positivity argument (Lemma~\ref{lem:halfdeg}) or by a recurrence (Section~\ref{sec:three}). The same reduction gives the general upper bounds.

For example, the case $k=0$ characterises the uniquely colourable graphs with $\AT(G)=\chi(G)$, which recovers a result of Hefetz \cite{Hefetz} (Corollary~\ref{cor:k0}). Theorem~\ref{thm:collapse} also gives upper bounds for graphs that are not uniquely colourable: adding suitable edges makes a graph uniquely colourable without changing $\chi$, and then $f(G)\le\AT(M_H)-\chi$ for the resulting graph $H$ (Corollary~\ref{cor:embed}). This is how Theorem~\ref{thm:intro-edges} is proved.

\subsection{Exact values}
Theorem~\ref{thm:collapse} turns the computation of $f$ into the computation of a coefficient of $P_{M_G}$, which for complete multipartite graphs is a constant term of a highly symmetric product. Table~\ref{tab:summary} collects the exact values obtained in this way. For complete multipartite graphs $K_{\vec m}=K_{m_1,\dots,m_\chi}$ the table refers to an explicit integer $D_0(\vec m)$, defined in Section~\ref{sec:multi}. Informally, $D_0(\vec m)+1$ is the lower bound for $\AT(M_G)$ obtained by counting the edges of $M_G$ inside each set of vertices, taking into account that certain pairs of vertices, called twins, must receive different exponents; it can be computed directly from $\vec m$ by an explicit formula. We call the resulting bound the \emph{twin-density bound}. The bipartite case deserves a separate statement. For a connected graph $C$ let $\beta(C)=e(C)-|V(C)|+1$ be its circuit rank.

\begin{theorem}\label{thm:bipartite}
For every bipartite graph $G$ with at least one edge,
\[
f(G)=\sum_C\Bigl\lfloor\frac{\beta(C)}{2}\Bigr\rfloor,
\]
where the sum is over the components $C$ of $G$.
\end{theorem}

\begin{table}[ht]
\centering\footnotesize
\renewcommand{\arraystretch}{1.22}
\begin{tabular}{@{}p{.26\linewidth}p{.36\linewidth}p{.14\linewidth}p{.17\linewidth}@{}}
\toprule
Graph $G$ & Value or bound for $f(G)$ & Result & Main tool\\
\midrule
uniquely $\chi$-colourable & $\AT(M_G)-\chi$ & Thm.~\ref{thm:collapse} & dual Cauchy identity\\
bipartite & $\sum_C\lfloor\beta(C)/2\rfloor$ & Thm.~\ref{thm:bipartite} & factorisation\\
$K_{m\star\chi}$ (equal parts) & $(\chi-1)\lfloor(m-1)^2/2\rfloor$ & Thm.~\ref{thm:equal} & Dyson identity\\
$K_{p,q,\dots,q}$ ($\chi-1$ parts of size $q$) & $D_1-\chi+1$ & Thm.~\ref{thm:morris} & Morris identity\\
$K_{m_1,m_2,m_3}$ (arbitrary parts) & $D_0-2\le f\le D_0-1$ & Cor.~\ref{cor:three} & three-term recurrence\\
$K_{m_1,m_2,m_3}$, $m_i\le60$ & $f=D_0-2$ & Cor.~\ref{cor:three}(b) & plus computation\\
$K_{2,2,m}$ & $\lceil(2m-4)/3\rceil$ & Cor.~\ref{cor:three} & log-concavity\\
\midrule
any $G$ on $n$ vertices & $\le n^2/8$, tight up to $O(n)$ & Cor.~\ref{cor:extremal} & torus positivity\\
$G(n,1/2)$ & $(\frac12+o(1))\,n\log_2n$ with high probability & Cor.~\ref{cor:random} & \\
\bottomrule
\end{tabular}
\caption{Main results. Here $D_0=D_0(\vec m)$ for complete multipartite graphs, as defined in Section~\ref{sec:multi}, and $D_1$ is given in Theorem~\ref{thm:morris}. The row for $m_i\le60$ is a finite computation (Section~\ref{sec:three}).}\label{tab:summary}
\end{table}

The lower bound in Theorem~\ref{thm:main} need not be attained. For example, it gives $1$ for two disjoint copies of $K_{2,3}$, whereas Theorem~\ref{thm:bipartite} gives $f(2K_{2,3})=2$. Here the two colourings of a component $C$ with $e(C)-|V(C)|$ odd cancel. For complete multipartite graphs, analogous cancellations arise from twins in $M_G$, and dense sub-multigraphs of $M_G$ give an additional obstruction; these are exactly the two effects measured by $D_0$. For instance, $f(K_{3,3,3})=4$ and $f(K_{2,5,8})=13$, while the lower bound in Theorem~\ref{thm:main} gives $3$ and $12$, respectively.

For complete tripartite graphs the result is sharper: Corollary~\ref{cor:three} places $f(K_{m_1,m_2,m_3})$ in the interval $[D_0(\vec m)-2,D_0(\vec m)-1]$, so the lower bound is within one of the true value for every choice of part sizes. Moreover, a computer check shows that the lower endpoint is attained whenever $m_1,m_2,m_3\le60$. For arbitrary complete multipartite graphs we conjecture that the twin and density obstructions account for the entire discrepancy (Conjecture~\ref{conj:M}); the conjecture is proved here when $\chi=2$, when all parts but one have the same size, and in the stated tripartite cases. Example~\ref{ex:3jzero} shows that the relevant coefficient can vanish for general uniquely $3$-colourable graphs, so the complete-multipartite hypothesis matters.

\subsection{Method}
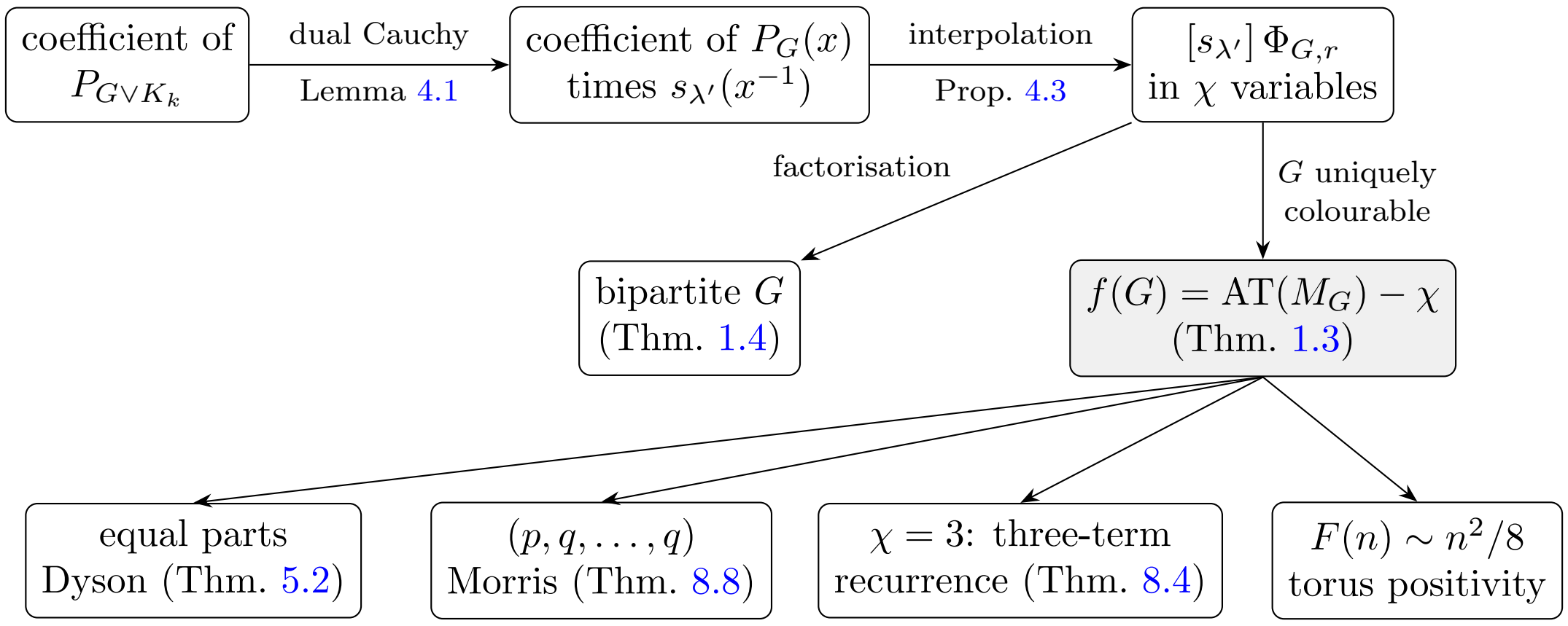
\begin{figure}[ht]
\centering
\begin{tikzpicture}[
  box/.style={draw,rounded corners=3pt,align=center,inner sep=4pt,font=\small,fill=#1},
  box/.default=white,
  lab/.style={font=\scriptsize,align=center},
  >={Stealth}]
\node[box] (c) at (0,0) {coefficient of\\ $P_{G\vee K_k}$};
\node[box] (p) at (5.1,0) {coefficient of $P_G(x)$\\ times $s_{\lambda'}(x^{-1})$};
\node[box] (s) at (10.3,0) {$[s_{\lambda'}]\,\Phi_{G,r}$\\ in $\chi$ variables};
\draw[->] (c) -- node[lab,above] {dual Cauchy} node[lab,below] {Lemma~\ref{lem:intY}} (p);
\draw[->] (p) -- node[lab,above] {interpolation} node[lab,below] {Prop.~\ref{prop:reduction}} (s);
\node[box=black!6] (m) at (10.3,-2.3) {$f(G)=\AT(M_G)-\chi$\\ (Thm.~\ref{thm:collapse})};
\draw[->] (s) -- node[lab,right] {$G$ uniquely\\ colourable} (m);
\node[box] (bip) at (5.1,-2.3) {bipartite $G$\\ (Thm.~\ref{thm:bipartite})};
\draw[->] (s.south west) -- node[lab,above left,pos=0.5] {factorisation} (bip.north east);
\node[box] (d) at (0.6,-4.5) {equal parts\\ Dyson (Thm.~\ref{thm:equal})};
\node[box] (mo) at (4.3,-4.5) {$(p,q,\dots,q)$\\ Morris (Thm.~\ref{thm:morris})};
\node[box] (t) at (8.1,-4.5) {$\chi=3$: three-term\\ recurrence (Thm.~\ref{thm:three})};
\node[box] (h) at (11.7,-4.5) {$F(n)\sim n^2/8$\\ torus positivity};
\foreach \x in {d,mo,t,h}{\draw[->] (m.south) -- (\x.north);}
\end{tikzpicture}
\caption{Structure of the proofs.}\label{fig:roadmap}
\end{figure}

The graph polynomial of $G\vee K_k$ factorises as $P_G(x)\cdot\prod_{v,i}(x_v-y_i)\cdot V(y)$, where $x$ are the variables of $G$ and $y$ those of the universal vertices. Figure~\ref{fig:roadmap} shows how the argument proceeds. Expanding the middle factor with the dual Cauchy identity for Schur functions separates $x$ from $y$, and extracting a coefficient in $y$ leaves a single term: every relevant coefficient of $P_{G\vee K_k}$ is a coefficient of $P_G(x)$ times one Schur polynomial in the reciprocals $1/x_v$. This Schur polynomial is indexed by a partition with at most $\chi$ parts, and Lagrange interpolation on $\chi$ points converts the coefficient into a Schur coefficient of a symmetric polynomial $\Phi_{G,r}$ in only $\chi$ variables. This polynomial is a weighted sum over the proper colourings of $G$ with $\chi$ colours; we call the weights, which are products of differences of colours, \emph{interpolation weights}. For uniquely colourable graphs there is essentially one such colouring, and the Schur coefficient becomes a constant term whose integrand is the graph polynomial of $M_G$. The colouring sums are those of the quantitative Combinatorial Nullstellensatz, and Appendix~\ref{app:qcn} gives a second proof of the reduction that applies the Nullstellensatz directly to $P_{G\vee K_k}$.

\subsection{Related work}
Hefetz \cite{Hefetz} expressed the Alon--Tarsi condition as a weighted sum over proper colourings; this is the special case of Theorem~\ref{thm:qcn} below with grids $\{0,1,\dots,a_v\}$. He used it to prove $\AT(G)=\chi(G)$ for uniquely colourable graphs satisfying certain edge-count and parity conditions, including $\AT(K_{2\star n})=n$. Kaul and Mudrock \cite{KM} applied the Combinatorial Nullstellensatz to DP-colouring, in particular to cones $K_1\vee G$ over certain bipartite graphs and over uniquely $3$-colourable graphs. Petrov and Gordeev \cite{PG} used almost central monomials of graph polynomials to bound the Alon--Tarsi numbers of Cartesian products. The Ohba number for choosability of complete bipartite graphs was investigated in \cite{CanoEtAl}, and paint and enumerative variants of Ohba's problem appear in \cite{CLMPTW,XYZ}. Several degeneracy-type parameters are known to bound the Alon--Tarsi number from above, for instance the variants of strict degeneracy studied by Zhou, Zhu and Zhu \cite{ZZZ}; we do not know whether any of them shows that $f(G)$ is finite. We are not aware of previous work on $\AT(G\vee K_k)$ for $k\ge2$, and the recent survey \cite{VSN} lists none.

\subsection{Organisation}
Section~\ref{sec:prelim} collects the tools. Section~\ref{sec:lower} proves the lower bounds, which come from twins and from edge density. Section~\ref{sec:reduction}, the technical core, eliminates the universal vertices. Section~\ref{sec:unique} proves Theorem~\ref{thm:collapse} and derives Theorems~\ref{thm:main} and \ref{thm:intro-edges} and the results on graphs with few vertices, and Section~\ref{sec:bip} treats bipartite graphs. Section~\ref{sec:multi} formulates Conjecture~\ref{conj:M} for complete multipartite graphs, and Section~\ref{sec:three} proves it in the cases listed above. Section~\ref{sec:open} collects open problems, and Appendix~\ref{app:qcn} gives a second proof of the main reduction via the Nullstellensatz. A reader interested mainly in the graph-theoretic results can read the statements in Sections~\ref{sec:lower} and \ref{sec:reduction}, then Sections~\ref{sec:unique}--\ref{sec:multi}; the algebra is confined to the proofs in Sections~\ref{sec:prelim}, \ref{sec:reduction} and \ref{sec:three}.

\section{Preliminaries}\label{sec:prelim}

We need three tools: Lagrange interpolation, in the form of Lemma~\ref{lem:dd} and the quantitative Combinatorial Nullstellensatz, which computes coefficients from values on grids; some basic facts about symmetric functions; and Dyson's constant term identity. Apart from Dyson's identity and the standard facts about symmetric functions quoted below, this section is self-contained.

\emph{Graphs.} We write $e(G)$ for the number of edges, $e(X,Y)$ for the number of edges between disjoint vertex sets $X$ and $Y$, $G[X]$ for the subgraph induced by $X$, $\alpha(G)$ for the independence number and $\Delta(G)$ for the maximum degree. We write $K_{m_1,\dots,m_\chi}$ for the complete multipartite graph with parts of sizes $m_1,\dots,m_\chi$, also denoted $K_{\vec m}$, and $K_{m\star\chi}$ for the case of $\chi$ parts of size $m$. We write $tH$ for $t$ disjoint copies of $H$ and $G\sqcup H$ for a disjoint union. A multigraph is \emph{loopless} if it has no loops; parallel edges are allowed.

\emph{Polynomials.} For a vector of variables $u=(u_1,u_2,\dots)$ and an integer vector $\nu$ with finitely many nonzero entries we write $u^\nu=\prod_iu_i^{\nu_i}$; for example $x^a=\prod_vx_v^{a_v}$ for $a\in\N^V$. We write $\coef{x^a}F$ for the coefficient of $x^a$ in $F$; $\one$ is the all-ones vector. For a Laurent polynomial $F$, $\CT F$ is its constant term, and $F(u^{-1})$ means $F(u_1^{-1},\dots,u_r^{-1})$. For a finite set $A\subseteq\Q$ and $t\in A$ we write $N_A(t)=\prod_{s\in A\setminus\{t\}}(t-s)$, and $[0,D]=\{0,1,\dots,D\}$.

\emph{Partitions.} A \emph{partition} $\lambda=(\lambda_1\ge\lambda_2\ge\dots\ge0)$ is a weakly decreasing sequence of nonnegative integers with finitely many nonzero terms. Its \emph{size} is $|\lambda|=\sum_i\lambda_i$ and its \emph{length} $\ell(\lambda)$ is the number of nonzero terms. The \emph{Young diagram} of $\lambda$ is the set of cells $\{(i,j):i\ge1,\ 1\le j\le\lambda_i\}$, drawn with $\lambda_i$ cells in row $i$; the cells $(1,j),(2,j),\dots$ in it form column $j$. The \emph{conjugate} partition $\lambda'$ is defined by $\lambda'_j=\#\{i:\lambda_i\ge j\}$, the length of column $j$; for example $(3,1)'=(2,1,1)$. Thus $\lambda'_1=\ell(\lambda)$ and $\ell(\lambda')=\lambda_1$. We write $(a^b)$ for the rectangular partition with $b$ parts equal to $a$, and $\lambda\subseteq\mu$ if $\lambda_i\le\mu_i$ for all $i$; thus $\lambda\subseteq(\chi^k)$ means that $\lambda$ has at most $k$ parts, each at most $\chi$. For partitions of the same size, $\nu$ is \emph{dominated} by $\mu$ if $\nu_1+\dots+\nu_i\le\mu_1+\dots+\mu_i$ for all $i$.

Let $h_r(S)=\sum_{i_1\le\dots\le i_r}s_{i_1}\cdots s_{i_r}$ denote the complete homogeneous symmetric polynomial of degree $r$ in the elements of a finite set $S=\{s_1,s_2,\dots\}$, with $h_0=1$ and $h_r=0$ for $r<0$.

\begin{lemma}\label{lem:dd}
If $S$ is a finite set of numbers and $j\ge0$, then $\sum_{t\in S}t^j/N_S(t)=h_{j-|S|+1}(S)$.
\end{lemma}

\begin{proof}
Compare coefficients of $z^j$ in the partial fraction expansion
\[
\frac{z^{|S|-1}}{\prod_{s\in S}(1-sz)}=\sum_{t\in S}\frac{1}{N_S(t)\,(1-tz)},
\]
an identity of rational functions in $z$ and the elements of $S$ (for $0\notin S$ it is the partial fraction decomposition, as both sides have simple poles at $z=1/t$ with equal residues and vanish at $z=\infty$; the general case follows by continuity in the elements of $S$); the left-hand side is $\sum_{r\ge0}h_r(S)z^{r+|S|-1}$.
\end{proof}

The following coefficient formula is due to Schauz \cite{SchauzQCN}, Laso\'n \cite{Lason} and Karasev and Petrov \cite{KP}; it refines the Combinatorial Nullstellensatz \cite{AlonCN}, which asserts only that $F$ cannot vanish on the whole grid when this coefficient is nonzero. The colouring sums of Section~\ref{sec:reduction} are sums of this form, Appendix~\ref{app:qcn} derives the main reduction directly from it, and the computations in Section~\ref{sec:multi} use it. Since its proof is three lines given Lemma~\ref{lem:dd}, we include it.

\begin{theorem}[Quantitative Combinatorial Nullstellensatz]\label{thm:qcn}
Let $F\in\Q[z_1,\dots,z_N]$ have degree at most $a_1+\dots+a_N$, and let $A_1,\dots,A_N\subseteq\Q$ with $|A_i|=a_i+1$. Then
\[
\coef{z^a}F=\sum_{z\in A_1\times\dots\times A_N}\frac{F(z)}{N_{A_1}(z_1)\cdots N_{A_N}(z_N)}.
\]
\end{theorem}

\begin{proof}
Both sides are linear in $F$, so we may take $F=z^c$ with $c_1+\dots+c_N\le a_1+\dots+a_N$. The right-hand side then factorises as $\prod_i\sum_{t\in A_i}t^{c_i}/N_{A_i}(t)=\prod_ih_{c_i-a_i}(A_i)$, by Lemma~\ref{lem:dd}. If $c=a$ this equals $1$. Otherwise $c_i<a_i$ for some $i$, because $\sum c_i\le\sum a_i$, and the product vanishes.
\end{proof}

\emph{Symmetric functions.} A polynomial in $r$ variables is \emph{symmetric} if it is unchanged by every permutation of the variables. A \emph{symmetric function} is a formal power series of bounded degree in infinitely many variables $u_1,u_2,\dots$ that is unchanged by every permutation of the variables; they form a ring $\Lambda$. Setting $u_{r+1}=u_{r+2}=\dots=0$ turns a symmetric function into a symmetric polynomial in $u_1,\dots,u_r$, its \emph{restriction} to $r$ variables. For any set or vector of variables $u$, $h_r(u)$ is as above, $e_r(u)=\sum_{i_1<\dots<i_r}u_{i_1}\cdots u_{i_r}$ and $p_r(u)=\sum_iu_i^r$, with $h_r=e_r=0$ for $r<0$. We write $\delta_r=(r-1,r-2,\dots,0)$, and $\Delta(u)=\prod_{i<j}(u_i-u_j)$ for the Vandermonde product of $u=(u_1,\dots,u_r)$; it should not be confused with the maximum degree $\Delta(G)$, which appears only in Section~\ref{sec:unique}. A \emph{semistandard tableau} of shape $\lambda$ is a filling of the cells of the Young diagram of $\lambda$ with positive integers that weakly increase along each row and strictly increase down each column. Its \emph{content} is the sequence $\nu=(\nu_1,\nu_2,\dots)$, where $\nu_i$ is the number of cells containing the entry $i$. For example, the filling of the shape $(2,1)$ with rows $1,1$ and $2$ has content $(2,1)$, while the fillings with rows $1,2$ and $3$, and with rows $1,3$ and $2$, both have content $(1,1,1)$. The \emph{Kostka number} $K_{\lambda\nu}$ is the number of semistandard tableaux of shape $\lambda$ and content $\nu$; thus $K_{(2,1),(1,1,1)}=2$. The \emph{Schur polynomial} in $r$ variables is
\[
s_\lambda(u_1,\dots,u_r)=\sum_Tu^{\nu(T)}=\sum_{\nu\in\N^r}K_{\lambda\nu}\,u^\nu,
\]
where $T$ ranges over the semistandard tableaux of shape $\lambda$ with entries in $\{1,\dots,r\}$ and $\nu(T)$ is the content of $T$. For example, $s_{(2,1)}(u_1,u_2)=u_1^2u_2+u_1u_2^2$, the two terms coming from the fillings with rows $1,1$ and $2$, and with rows $1,2$ and $2$. The same sum over all semistandard tableaux of shape $\lambda$, in the variables $u_1,u_2,\dots$, is the \emph{Schur function} $s_\lambda$; its restriction to $r$ variables is the Schur polynomial. We use the following standard facts without proof; all can be found in Chapter~7 of \cite{EC2} or Chapter~I of \cite{Mac}.
\begin{enumerate}[label=\textup{(S\arabic*)}]
\item Each $s_\lambda$ is a symmetric function, and the $s_\lambda$, over all partitions $\lambda$, form a basis of $\Lambda$. Moreover $K_{\lambda\lambda}=1$, and for a partition $\nu$ we have $K_{\lambda\nu}=0$ unless $\nu$ is dominated by $\lambda$.
\item The restriction of $s_\lambda$ to $r$ variables is $0$ if $\ell(\lambda)>r$, and the restrictions of the $s_\lambda$ with $\ell(\lambda)\le r$ form a basis of the symmetric polynomials in $r$ variables.
\item \emph{(Bialternant formula.)} For a vector $\alpha\in\N^r$ let $a_\alpha=\det(u_i^{\alpha_j})_{i,j=1}^r$. Then $a_{\delta_r}=\Delta(u)$, and $s_\lambda(u_1,\dots,u_r)=a_{\lambda+\delta_r}/a_{\delta_r}$ whenever $\ell(\lambda)\le r$.
\item \emph{(Young's rule.)} For every sequence $\beta=(\beta_1,\dots,\beta_N)$ of integers, $h_{\beta_1}\cdots h_{\beta_N}=\sum_\nu K_{\nu\beta}\,s_\nu$, where $K_{\nu\beta}$ is the number of semistandard tableaux of shape $\nu$ and content $\beta$ (and $K_{\nu\beta}=0$ if some $\beta_i<0$).
\item \emph{(Dual Cauchy identity.)} For variables $u_1,\dots,u_r$ and $w_1,\dots,w_N$,
\[
\prod_{i=1}^r\prod_{j=1}^N(1+u_iw_j)=\sum_\mu s_\mu(u_1,\dots,u_r)\,s_{\mu'}(w_1,\dots,w_N),
\]
the sum over all partitions $\mu$; by (S2), only finitely many terms are nonzero.
\end{enumerate}
For a symmetric polynomial $\varphi$ in $r$ variables and a partition $\lambda$ with at most $r$ parts, $[s_\lambda]\varphi$ denotes the coefficient of $s_\lambda$ when $\varphi$ is expanded in the basis of (S2). It is a number; the variables only serve to write $\varphi$ down.

\begin{lemma}\label{lem:schur}
Let $\varphi$ be a symmetric polynomial in $u=(u_1,\dots,u_r)$ and $\lambda$ a partition with at most $r$ parts. Then
\[
[s_\lambda]\varphi=\coef{u^{\lambda+\delta_r}}\bigl(\Delta(u)\varphi(u)\bigr)=\frac1{r!}\CT\Bigl[\varphi(u)\,s_\lambda(u^{-1})\,\Delta(u)\Delta(u^{-1})\Bigr],
\]
where $\CT$ denotes the constant term, and $\Delta(u)\Delta(u^{-1})=(-1)^{\binom r2}\Delta(u)^2\prod_iu_i^{-(r-1)}=\prod_{i\ne j}(1-u_i/u_j)$.
\end{lemma}

\begin{proof}
Write $c_\mu=[s_\mu]\varphi$, so that $\varphi=\sum_\mu c_\mu s_\mu$, the sum over $\mu$ with $\ell(\mu)\le r$ by (S2). By (S3), $\Delta\varphi=\sum_\mu c_\mu a_{\mu+\delta_r}$. Expanding the determinant,
\[
a_\alpha(u)=\sum_{\sigma}\sgn(\sigma)\,u_1^{\alpha_{\sigma(1)}}\cdots u_r^{\alpha_{\sigma(r)}},
\]
the sum over all permutations $\sigma$ of $\{1,\dots,r\}$. If $\alpha$ is strictly decreasing, its entries are distinct, and the only term whose exponent vector is strictly decreasing is the term $u^\alpha$ with $\sigma=\mathrm{id}$. Hence the monomial $u^{\lambda+\delta_r}$ occurs in $a_{\mu+\delta_r}$ only when $\mu=\lambda$, with coefficient $1$; this proves the first equality.

For the second, apply (S3) to the variables $u^{-1}$: $s_\lambda(u^{-1})\Delta(u^{-1})=a_{\lambda+\delta_r}(u^{-1})$. Together with the expansion of $\Delta\varphi$ this gives
\[
\varphi(u)\,s_\lambda(u^{-1})\,\Delta(u)\Delta(u^{-1})=\sum_\mu c_\mu\,a_{\mu+\delta_r}(u)\,a_{\lambda+\delta_r}(u^{-1}),
\]
so it suffices to show that $\CT[a_\alpha(u)a_\beta(u^{-1})]=r!\,[\alpha=\beta]$ for strictly decreasing $\alpha,\beta$. Expanding both determinants, $a_\alpha(u)a_\beta(u^{-1})=\sum_{\sigma,\tau}\sgn(\sigma\tau)\prod_iu_i^{\alpha_{\sigma(i)}-\beta_{\tau(i)}}$. A term is constant only if $\alpha_{\sigma(i)}=\beta_{\tau(i)}$ for all $i$, which for strictly decreasing vectors forces $\alpha=\beta$ and $\sigma=\tau$; there are then $r!$ such terms, each equal to $1$. Finally, $(u_i-u_j)(u_i^{-1}-u_j^{-1})=-(u_i-u_j)^2/(u_iu_j)=(1-u_i/u_j)(1-u_j/u_i)$ for each pair $i<j$, and taking the product over all pairs gives the two expressions for $\Delta(u)\Delta(u^{-1})$.
\end{proof}

The last tool is a classical evaluation of a constant term. It was conjectured by Dyson \cite{Dyson} and proved in 1962 by Wilson \cite{Wilson} and, independently, by Gunson, whose proof was not published (see \cite{FW}); Good \cite{Good} gave a one-page proof by Lagrange interpolation, and Zeilberger \cite{Zeil} a combinatorial proof by counting tournaments.

\begin{theorem}[Dyson's constant term identity]\label{thm:dyson}
For nonnegative integers $a_1,\dots,a_r$,
\[
\CT\prod_{i\ne j}\Bigl(1-\frac{u_i}{u_j}\Bigr)^{a_i}=\frac{(a_1+\dots+a_r)!}{a_1!\cdots a_r!}.
\]
\end{theorem}

\section{Lower bounds: twins and density}\label{sec:lower}

This section sets up notation for $G\vee K_k$ and proves the lower bounds. There are two obstructions to a small Alon--Tarsi number: there may be too many edges for the available exponents (Proposition~\ref{prop:density}), and twin vertices must receive distinct exponents (Lemma~\ref{lem:twin}).

In $G\vee K_k$ the $k$ universal vertices are placed after $V(G)$ in the vertex order, and their variables are denoted $y_1,\dots,y_k$. Thus
\[
P_{G\vee K_k}=P_G(x)\cdot\prod_{v\in V(G)}\prod_{i=1}^{k}(x_v-y_i)\cdot V(y),\qquad V(y)=\prod_{1\le i<j\le k}(y_i-y_j).
\]
In Sections~\ref{sec:lower}--\ref{sec:bip} and in Appendix~\ref{app:qcn} we write $D=\chi+k-1$; the aim is to find a nonzero monomial of $P_{G\vee K_k}$ in which all exponents are at most $D$.

\begin{lemma}\label{lem:gap}
$\AT(G\vee K_1)\le\AT(G)+1$. Hence $\AT(G\vee K_k)=\chi+k$ for every $k\ge f(G)$. Moreover $f(G)\le f(H)$ whenever $G\subseteq H$ and $\chi(G)=\chi(H)$.
\end{lemma}

\begin{proof}
If $\coef{x^a}P_G\ne0$, then the coefficient of $x^{a+\one}y^0$ in $P_{G\vee K_1}=P_G(x)\prod_v(x_v-y)$ equals $\coef{x^a}P_G$, since $y^0$ forces the term $x_v$ from every factor; its exponents are at most $\AT(G)$. Applied to $G\vee K_k$ in place of $G$, this gives $\AT(G\vee K_{k+1})\le\AT(G\vee K_k)+1$. Since always $\AT(G\vee K_k)\ge\chi(G\vee K_k)=\chi+k$, the gap $\AT(G\vee K_k)-\chi-k$ is nonnegative and does not increase with $k$, which proves the second statement. The last statement follows because $\AT$ is monotone under taking subgraphs: if $G\subseteq H$, then $P_H$ is $\pm P_G$ times a product of linear factors $x_u-x_w$ (reordering the vertices only changes the sign), and if every monomial of a polynomial $F$ has an exponent at least $d$, then so does every monomial of $(x_u-x_w)F$.
\end{proof}

\begin{lemma}\label{lem:twin}
Let $u$ and $u'$ be adjacent vertices of a graph $H$ with the same closed neighbourhood. Then interchanging $x_u$ and $x_{u'}$ changes the sign of $P_H$. Consequently, interchanging the exponents of $x_u$ and $x_{u'}$ in a monomial changes the sign of its coefficient, and monomials in which these exponents are equal have coefficient $0$.
\end{lemma}

\begin{proof}
Every other vertex $w$ is adjacent to both $u$ and $u'$ or to neither. In the first case the product of the two factors of $P_H$ corresponding to $wu$ and $wu'$ is invariant under the interchange, whatever the relative order of $w,u,u'$. The factor corresponding to $uu'$ changes sign.
\end{proof}

We call such vertices \emph{twins}. The universal vertices of $G\vee K_k$ are pairwise twins, so their exponents in a nonzero monomial are distinct, and after permuting the $y_i$ (which changes only the sign) every nonzero monomial of $P_{G\vee K_k}$ has the form $x^\alpha y^b$ with $b_1>b_2>\dots>b_k$. We write $b=\delta_k+\lambda$, where $\lambda$ is a partition with at most $k$ parts. If all exponents are at most $D$, then $\lambda\subseteq(\chi^k)$. Accordingly, we call a monomial $x^\alpha y^{\delta_k+\lambda}$ \emph{admissible} if $\lambda\subseteq(\chi^k)$, all $\alpha_v\le D$, and its total degree equals $\deg P_{G\vee K_k}=e(G)+nk+\binom k2$; its coefficient in $P_{G\vee K_k}$ may be zero. By the above, $\AT(G\vee K_k)\le\chi+k$ if and only if some admissible monomial has a nonzero coefficient. Comparing degrees shows that every admissible monomial satisfies
\begin{equation}\label{eq:budget}
\sum_{v\in V(G)}(D-\alpha_v)+\bigl(k\chi-|\lambda|\bigr)=(\chi-1)n+k\chi-e(G).
\end{equation}
The quantities $r_v=D-\alpha_v\ge0$ are the \emph{deficits} of the monomial, and we call the right-hand side $(\chi-1)n+k\chi-e(G)$ of \eqref{eq:budget} the \emph{budget}: the deficits and the shortfall $k\chi-|\lambda|$ must add up to it.

\begin{proposition}\label{prop:density}
If $\AT(G\vee K_k)=\chi+k$, then $e(G')\le(\chi-1)|V(G')|+k\chi$ for every subgraph $G'$ of $G$.
\end{proposition}

\begin{proof}
By monotonicity $P_{G'\vee K_k}$ has an admissible nonzero monomial, and the left-hand side of \eqref{eq:budget}, written for $G'$, is nonnegative.
\end{proof}

The density bound uses only the chromatic number of $G$. A more refined bound compares $G$ with smaller graphs of smaller chromatic number.

\begin{proposition}\label{prop:subjoin}
If $G$ contains $G'\vee K_s$ as a subgraph and $\chi(G)=\chi(G')+s$, then $f(G)\ge f(G')-s$.
\end{proposition}

\begin{proof}
Let $k=f(G)$. Then $G\vee K_k$ contains $G'\vee K_{k+s}$, so $\AT(G'\vee K_{k+s})\le\AT(G\vee K_k)=\chi(G)+k=\chi(G')+k+s$. Hence $f(G')\le k+s$.
\end{proof}

\begin{example}\label{ex:gaps}
The density bound of Proposition~\ref{prop:density} is not within an additive constant of $f$; part (a) shows this already for $\chi=2$.
\begin{enumerate}[label=\textup{(\alph*)}]
\item For $t$ disjoint copies of $K_{2,3}$, Theorem~\ref{thm:bipartite} gives $f=t$, while the density bound is $t/2$.
\item The graph $K_{1,m,m}$ contains $K_{m,m}\vee K_1$, so Proposition~\ref{prop:subjoin} gives $f(K_{1,m,m})\ge\lfloor(m-1)^2/2\rfloor-1$. The density bound gives only $(m^2-2m-2)/3$, attained at $G'=G$. In fact $f(K_{1,m,m})=\lfloor(m-1)^2/2\rfloor-1$ for $m\ge3$, by Theorem~\ref{thm:morris} below.
\end{enumerate}
Similarly, $f(K_{2,5,8})\ge f(K_{5,8})-1=13$, which exceeds the density bound $12$.
\end{example}

\section{Eliminating the universal vertices}\label{sec:reduction}

We now show that every admissible coefficient of $P_{G\vee K_k}$ is a Schur coefficient of a symmetric polynomial in only $\chi$ variables, which depends on $G$ and on the deficits of the vertices of $G$ (Proposition~\ref{prop:reduction}). The number $k$ remains only in the shape of the partition indexing that Schur coefficient. There are two steps. First, the factors of $P_{G\vee K_k}$ that involve the universal vertices are expanded with the dual Cauchy identity (S5); this removes the variables $y$ (Lemma~\ref{lem:intY}). Second, Lagrange interpolation on $\chi$ points turns the result into a Schur coefficient in $\chi$ variables (Lemma~\ref{lem:Phipoly}).

Throughout this section we fix $k$ and an admissible monomial $x^\alpha y^{\delta_k+\lambda}$ with deficits $r_v=D-\alpha_v$, and we write $x^{-1}=(x_1^{-1},\dots,x_n^{-1})$.

\begin{lemma}\label{lem:intY}
We have
\[
\coef{x^\alpha y^{\delta_k+\lambda}}P_{G\vee K_k}=(-1)^{|\lambda|}\,\CT\Bigl[P_G(x)\,s_{\lambda'}(x^{-1})\prod_vx_v^{\,r_v-\chi+1}\Bigr].
\]
\end{lemma}

\begin{proof}
We compute in the ring of polynomials in $y$ whose coefficients are Laurent polynomials in $x$; it contains $\Q[x,y]$, and taking the coefficient of a monomial in $y$ commutes with this inclusion. Factor out $x_v$ from each $x_v-y_i$, and substitute $u_i=y_i$ and $w_v=-1/x_v$ in the identity (S5):
\begin{align*}
\prod_v\prod_i(x_v-y_i)&=(x_1\cdots x_n)^k\prod_{i,v}\Bigl(1-\frac{y_i}{x_v}\Bigr)\\
&=(x_1\cdots x_n)^k\sum_\mu(-1)^{|\mu|}s_\mu(y)\,s_{\mu'}(x^{-1}).
\end{align*}
Here we used that $s_{\mu'}$ is homogeneous of degree $|\mu|$, and by (S2) only partitions $\mu$ with at most $k$ parts contribute. Multiplying by $P_G(x)$ and $V(y)=\Delta(y)$ gives
\[
P_{G\vee K_k}=\sum_\mu(-1)^{|\mu|}\,P_G(x)\,(x_1\cdots x_n)^k\,s_{\mu'}(x^{-1})\cdot\Delta(y)\,s_\mu(y).
\]
By Lemma~\ref{lem:schur} with $r=k$ and $\varphi=s_\mu$, the coefficient of $y^{\delta_k+\lambda}$ in $\Delta(y)s_\mu(y)$ is $[s_\lambda]s_\mu$, which is $1$ if $\mu=\lambda$ and $0$ otherwise. Hence the coefficient of $y^{\delta_k+\lambda}$ in $P_{G\vee K_k}$ is
\[
(-1)^{|\lambda|}P_G(x)\,(x_1\cdots x_n)^k\,s_{\lambda'}(x^{-1}).
\]
This is a polynomial in $x$, as it must be: every monomial of $s_{\lambda'}(x^{-1})$ is $x^{-\nu}$ with $\nu$ the content of a semistandard tableau of shape $\lambda'$, each value occurs at most once in each of the $\lambda'_1=\ell(\lambda)\le k$ columns of such a tableau, and so $\nu_v\le k$. Its coefficient of $x^\alpha$ is the constant term of this polynomial times $x^{-\alpha}$, and $k-\alpha_v=k-D+r_v=r_v-\chi+1$.
\end{proof}

Thus the $k$ universal vertices are absorbed into the single Laurent polynomial $s_{\lambda'}(x^{-1})$ in the variables of $G$. The partition $\lambda'$ has $\lambda_1\le\chi$ parts, so it indexes a Schur polynomial in $\chi$ variables; the next lemma moves the computation to $\chi$ variables.

\emph{The colouring polynomial.} For a set $S$ of $\chi$ numbers put
\begin{equation}\label{eq:Phi}
\Phi_{G,r}(S)=\sum_{c\colon V(G)\to S}\frac{P_G(c)\prod_vc_v^{\,r_v}}{\prod_vN_S(c_v)}.
\end{equation}
Since $P_G(c)=0$ unless $c$ is a proper colouring, the sum runs in effect over the proper colourings of $G$ with colours in $S$; as $|S|=\chi=\chi(G)$, each of them uses every colour of $S$. The sum is the right-hand side of Theorem~\ref{thm:qcn} for the polynomial $P_G(x)\prod_vx_v^{r_v}$ with every grid equal to $S$. The degree condition of Theorem~\ref{thm:qcn} fails here, so the sum is not a single coefficient; instead it is a symmetric polynomial in the elements of $S$.

\begin{lemma}\label{lem:Phipoly}
There is a symmetric polynomial $\Phi_{G,r}(u_1,\dots,u_\chi)$ whose value at the elements of every $\chi$-element set $S$, taken in any order, is $\Phi_{G,r}(S)$. For every partition $\nu$ with at most $\chi$ parts, its Schur coefficient is
\[
[s_\nu]\,\Phi_{G,r}=\CT\Bigl[P_G(x)\,s_\nu(x^{-1})\prod_vx_v^{\,r_v-\chi+1}\Bigr].
\]
\end{lemma}

\begin{proof}
Write $P_G(x)\prod_vx_v^{r_v}=\sum_af_ax^a$. For each monomial $x^a$ the sum over $c\colon V(G)\to S$ splits into one sum per vertex, so by Lemma~\ref{lem:dd} with $|S|=\chi$,
\[
\Phi_{G,r}(S)=\sum_af_a\prod_v\Bigl(\sum_{t\in S}\frac{t^{a_v}}{N_S(t)}\Bigr)=\sum_af_a\prod_vh_{a_v-\chi+1}(S).
\]
So $\Phi_{G,r}(u)=\sum_af_a\prod_vh_{a_v-\chi+1}(u)$ has the required values. It is the only polynomial with these values, because a polynomial that vanishes at every point with distinct coordinates is zero; this is what allows us to compute $\Phi_{G,r}(u)$ from its values, as in \eqref{eq:PhiU} below. Put $\beta=a-(\chi-1)\one$. By (S4) and (S2), the coefficient of $s_\nu$ in $\prod_vh_{\beta_v}(u)$ is the Kostka number $K_{\nu\beta}$, and by the definition of $s_\nu$ this is the coefficient of $x^\beta$ in $s_\nu(x_1,\dots,x_n)$, that is, $\CT[x^a\,s_\nu(x^{-1})\prod_vx_v^{-\chi+1}]$; both are $0$ if some $\beta_v<0$. Summing over $a$ with weights $f_a$ proves the formula.
\end{proof}

\begin{proposition}\label{prop:reduction}
For every admissible monomial $x^\alpha y^{\delta_k+\lambda}$ with deficits $r$,
\[
\coef{x^\alpha y^{\delta_k+\lambda}}P_{G\vee K_k}=(-1)^{|\lambda|}\,[s_{\lambda'}]\,\Phi_{G,r}.
\]
\end{proposition}

\begin{proof}
Since $\lambda\subseteq(\chi^k)$, the partition $\lambda'$ has at most $\chi$ parts. Apply Lemma~\ref{lem:Phipoly} with $\nu=\lambda'$ and compare with Lemma~\ref{lem:intY}.
\end{proof}

So the coefficient depends on $G$ only through the colouring polynomial $\Phi_{G,r}$, and the $k$ universal vertices are replaced by a Schur coefficient in $\chi$ variables.

\begin{remark}[The reduction via the Nullstellensatz]\label{rem:qcn}
Proposition~\ref{prop:reduction} can also be proved by applying Theorem~\ref{thm:qcn} to $P_{G\vee K_k}$ itself, and this explains its shape. Regard $[0,D]$ as a palette of $\chi+k$ colours. In the resulting grid sum, the universal vertices receive a set $Y$ of $k$ distinct colours, the factors $x_v-y_i$ cancel the part of the interpolation weights that involves $Y$, and $G$ is coloured from the complementary $\chi$ colours. Passing from this sum over $k$-sets to a Schur coefficient in $\chi$ variables uses the involution $\omega$ on symmetric functions; the dual Cauchy identity in Lemma~\ref{lem:intY} does this in one step. The details are in Appendix~\ref{app:qcn}.
\end{remark}

\section{The collapse theorem and its consequences}\label{sec:unique}

For uniquely colourable graphs, the colouring polynomial from Section~\ref{sec:reduction} is a single symmetrised product. This makes the reduction explicit and proves Theorem~\ref{thm:collapse}. We then apply the collapse theorem in two ways: Dyson's identity gives the exact answer for equal part sizes, while a torus-positivity argument gives a bound in terms of the edges leaving each colour class.

Let $G$ be uniquely $\chi$-colourable with colour classes $V_1,\dots,V_\chi$, $m_i=|V_i|$ and $e_{ij}=e(V_i,V_j)$. For deficits $r$ put $\rho_i=\sum_{v\in V_i}r_v$, the total deficit of the class $V_i$. The proper colourings of $G$ onto $\{u_1,\dots,u_\chi\}$ are $c_v=u_{\sigma(i)}$ for $v\in V_i$, where $\sigma$ ranges over the permutations of $\{1,\dots,\chi\}$. Writing $\xi_i=u_{\sigma(i)}$ for the colour of $V_i$ and using
\[
\prod_{i}\prod_{j\ne i}(\xi_i-\xi_j)^{m_i}=\prod_{i<j}(-1)^{m_j}(\xi_i-\xi_j)^{m_i+m_j}
\]
for the interpolation weights, one finds
\begin{equation}\label{eq:PhiU}
\Phi_{G,r}(u)=\varepsilon\sum_{\sigma}g(u_{\sigma(1)},\dots,u_{\sigma(\chi)}),\qquad g(\xi)=\prod_{i<j}(\xi_i-\xi_j)^{e_{ij}-m_i-m_j}\prod_i\xi_i^{\rho_i},
\end{equation}
with a sign $\varepsilon$ depending only on $G$ and the vertex order. Note that $e_{ij}-m_i-m_j=\mu_{ij}-2\ge-1$. Informally, each colour class $V_i$ behaves like a single vertex $u_i$, and the edges of $G$ between $V_i$ and $V_j$, minus the interpolation weights, leave $\mu_{ij}-2$ factors $u_i-u_j$; the factor $\Delta(u)\Delta(u^{-1})$ in Lemma~\ref{lem:schur} supplies the missing two.

\begin{proof}[Proof of Theorem~\ref{thm:collapse}]
Let $x^\alpha y^{\delta_k+\lambda}$ be admissible. By Proposition~\ref{prop:reduction}, Lemma~\ref{lem:schur} and \eqref{eq:PhiU}, its coefficient is, up to sign,
\[
\frac1{\chi!}\CT\Bigl[\Phi_{G,r}(u)\,s_{\lambda'}(u^{-1})\,\Delta(u)\Delta(u^{-1})\Bigr]=\CT\Bigl[g(u)\,s_{\lambda'}(u^{-1})\,\Delta(u)\Delta(u^{-1})\Bigr],
\]
since the weight is symmetric and each product $g(\sigma u)\Delta(u)\Delta(u^{-1})$ is a Laurent polynomial, so that its constant term does not change when the variables are permuted. (Every exponent $\mu_{ij}-2$ in $g$ is at least $-1$, and $\Delta(u)\Delta(u^{-1})$ is $(u_i-u_j)^2$ times a Laurent polynomial.) By Lemma~\ref{lem:schur}, $g(u)\Delta(u)\Delta(u^{-1})=\pm P_{M_G}(u)\,u^\rho\prod_iu_i^{-(\chi-1)}$, and therefore
\begin{equation}\label{eq:expansion}
\coef{x^\alpha y^{\delta_k+\lambda}}P_{G\vee K_k}=\pm\sum_\nu K_{\lambda'\nu}\,\coef{u^{\nu+(\chi-1)\one-\rho}}P_{M_G}(u).
\end{equation}
This uses $s_{\lambda'}(u^{-1})=\sum_\nu K_{\lambda'\nu}u^{-\nu}$, where $\nu$ ranges over the contents of semistandard tableaux of shape $\lambda'$. Each $i$ occurs at most once in each column of such a tableau, so $\nu_i\le\lambda'_1\le k$, and every monomial on the right has all exponents at most $D$.

Suppose first that $\AT(M_G)>\chi+k$. Then all coefficients on the right of \eqref{eq:expansion} vanish, hence every admissible coefficient of $P_{G\vee K_k}$ vanishes and $\AT(G\vee K_k)>\chi+k$.

Conversely, let $\coef{u^\gamma}P_{M_G}\ne0$ with $\gamma_i\le D$ for all $i$. Take $\lambda=(\chi^k)$, so that $s_{\lambda'}(u^{-1})=(u_1\cdots u_\chi)^{-k}$ and \eqref{eq:expansion} has a single term, and put $\rho_i=D-\gamma_i$. As $e(M_G)=e(G)-(\chi-1)n+\chi(\chi-1)$, we get $\sum_i\rho_i=\chi D-e(M_G)=(\chi-1)n+k\chi-e(G)$, which is the budget \eqref{eq:budget} for $|\lambda|=k\chi$. Distribute the deficit $\rho_i\le D$ arbitrarily over the vertices of $V_i$. The resulting monomial is admissible, and by \eqref{eq:expansion} its coefficient is $\pm\coef{u^\gamma}P_{M_G}\ne0$.
\end{proof}

\begin{corollary}\label{cor:k0}
A uniquely $\chi$-colourable graph $G$ satisfies $\AT(G)=\chi$ if and only if $P_{M_G}$ has a nonzero monomial with all exponents at most $\chi-1$. In particular, $\AT(G)=\chi$ whenever $e(G)=(\chi-1)n-\binom\chi2$.
\end{corollary}

\begin{proof}
The first statement is the case $k=0$ of Theorem~\ref{thm:collapse}. If $e(G)=(\chi-1)n-\binom\chi2$, then $\sum_{i<j}\mu_{ij}=\binom\chi2$; since every $\mu_{ij}\ge1$, $M_G$ is the simple graph $K_\chi$. Its graph polynomial is $\Delta(u)$, in which $u^{\delta_\chi}$ has coefficient $1$, so $\AT(K_\chi)\le\chi$, and $\AT(K_\chi)\ge\chi(K_\chi)=\chi$.
\end{proof}

The second statement of Corollary~\ref{cor:k0} is Theorem~1.8 of \cite{Hefetz}.

\begin{theorem}\label{thm:equal}
For all $m\ge1$ and $\chi\ge2$,
\[
f(K_{m\star\chi})=(\chi-1)\Bigl\lfloor\frac{(m-1)^2}{2}\Bigr\rfloor .
\]
\end{theorem}

\begin{proof}
Let $b=\lfloor(m-1)^2/2\rfloor$, $k=(\chi-1)b$ and $D=\chi+k-1$. All multiplicities of $M=M_{K_{m\star\chi}}$ equal $\mu=(m-1)^2+1$, so $P_M=\Delta(u)^\mu$. We use repeatedly that, by Lemma~\ref{lem:schur},
\begin{equation}\label{eq:dysonform}
\Delta(u)^{2c}\,(u_1\cdots u_\chi)^{-(\chi-1)c}=\pm\prod_{i\ne j}\Bigl(1-\frac{u_i}{u_j}\Bigr)^{c}\qquad(c\ge0).
\end{equation}

\emph{$m$ even.} Then $\mu=2(b+1)$ and $D=(\chi-1)(b+1)$. The lower bound $f\ge k$ is Proposition~\ref{prop:density}. By \eqref{eq:dysonform} and Theorem~\ref{thm:dyson}, $\coef{u_1^D\cdots u_\chi^D}P_M=\pm(\chi(b+1))!/((b+1)!)^\chi\ne0$, so $\AT(M)\le D+1=\chi+k$, and Theorem~\ref{thm:collapse} gives $f\le k$.

\emph{$m$ odd.} Then $\mu=2b+1$ is odd, so interchanging two variables multiplies $P_M=\Delta(u)^\mu$ by $(-1)^\mu=-1$. As in Lemma~\ref{lem:twin}, in every nonzero monomial of $P_M$ the $\chi$ exponents are therefore distinct. If they are at most $D'$, their sum $\binom\chi2(2b+1)$ is at most $\chi D'-\binom\chi2$, whence $D'\ge(\chi-1)(b+1)=D$. Thus $\AT(M)\ge\chi+k$, and $f\ge k$ by Theorem~\ref{thm:collapse}. Conversely, let $\gamma=\delta_\chi+(k^\chi)$, whose largest entry is $D$. By Lemma~\ref{lem:schur}, \eqref{eq:dysonform} and Theorem~\ref{thm:dyson},
\[
\coef{u^\gamma}\Delta^{2b+1}=[s_{(k^\chi)}]\Delta^{2b}=\frac1{\chi!}\CT\Bigl[\Delta^{2b}(u_1\cdots u_\chi)^{-k}\Delta(u)\Delta(u^{-1})\Bigr]=\pm\frac{1}{\chi!}\,\frac{(\chi(b+1))!}{((b+1)!)^\chi}\ne0 .
\]
Hence $\AT(M)\le\chi+k$ and $f\le k$.
\end{proof}

For $m=2$ Theorem~\ref{thm:equal} gives $\AT(K_{2\star\chi})=\chi$, which is Corollary~1.10 of \cite{Hefetz}, and for $m=3$ it gives $f(K_{3\star\chi})=2(\chi-1)$.

\begin{proof}[Proof of Theorem~\ref{thm:main}]
If $\chi=1$, then $G$ is edgeless and $f(G)=0$, so the assertion is immediate. Assume henceforth that $\chi\ge2$. The lower bound is Proposition~\ref{prop:density}. For the upper bound, enlarge each colour class to size $m$; this embeds $G$ into $K_{m\star\chi}$ without changing the chromatic number. By Lemma~\ref{lem:gap} and Theorem~\ref{thm:equal}, $f(G)\le f(K_{m\star\chi})=(\chi-1)\lfloor(m-1)^2/2\rfloor$. For even $m$ this equals $(\chi-1)m(m-2)/2$, which is the lower bound for $G=K_{m\star\chi}$.
\end{proof}

More generally, Theorem~\ref{thm:collapse} bounds $f$ for every graph, by embedding it into a uniquely colourable graph with the same chromatic number.

\begin{corollary}\label{cor:embed}
Let $G$ be a graph with $\chi(G)=\chi$.
\begin{enumerate}[label=\textup{(\alph*)}]
\item If $H\supseteq G$ is uniquely $\chi$-colourable and $\chi(H)=\chi$, then $f(G)\le f(H)=\AT(M_H)-\chi$.
\item Let $V_1,\dots,V_\chi$ be the classes of a proper $\chi$-colouring of $G$, with $m_i=|V_i|$ and $e_{ij}=e(V_i,V_j)$. Choose $z_i\in V_i$ for each $i$, and let $G'$ be obtained from $G$ by joining each $z_i$ to every vertex outside $V_i$. Then $G'$ is uniquely $\chi$-colourable with colour classes $V_1,\dots,V_\chi$, and
\[
f(G)\le\AT(M_{G'})-\chi,
\]
where $M_{G'}$ has multiplicities $\mu'_{ij}=e(G'[V_i\cup V_j])-m_i-m_j+2$ satisfying $1\le\mu'_{ij}\le e_{ij}+1$.
\end{enumerate}
In particular, $f(G)\le\AT(M_{K_{\vec m}})-\chi$ for the vector $\vec m$ of class sizes of any proper $\chi$-colouring of $G$.
\end{corollary}

\begin{proof}
(a) follows from Lemma~\ref{lem:gap} and Theorem~\ref{thm:collapse}. For (b), the vertices $z_1,\dots,z_\chi$ of $G'$ form a clique, and every $v\in V_i$ is adjacent to all $z_j$ with $j\ne i$. Hence every proper $\chi$-colouring of $G'$ gives $v$ the colour of $z_i$, so $G'$ is uniquely $\chi$-colourable with the classes $V_i$, and $\chi(G')=\chi$ because $G'$ contains the clique $\{z_1,\dots,z_\chi\}$. At most $m_i+m_j-1$ edges were added between $V_i$ and $V_j$, which gives $\mu'_{ij}\le e_{ij}+1$, and $\mu'_{ij}\ge1$ holds for every uniquely colourable graph. Now apply (a). The last statement is (a) with $H=K_{\vec m}$, which contains $G$ and is uniquely $\chi$-colourable.
\end{proof}

The graph $G'$ of (b) adds few edges, so $M_{G'}$ records the edges of $G$ between the colour classes rather than only the class sizes; Theorem~\ref{thm:edges} below exploits this.

\subsection*{A bound via the edges between colour classes}
Theorem~\ref{thm:main} is efficient when the colour classes are balanced. The next bound uses only the numbers of edges leaving each colour class. It rests on a positivity argument for multigraphs.

\begin{lemma}\label{lem:halfdeg}
Let $M$ be a multigraph on $\{1,\dots,r\}$ with $\mu_{ij}$ edges between $i$ and $j$. Then
\[
\AT(M)\le\max_i\sum_{j\ne i}\Bigl\lceil\frac{\mu_{ij}}2\Bigr\rceil+1 .
\]
\end{lemma}

\begin{proof}
Let $c_{ij}=\lceil\mu_{ij}/2\rceil$ and $w_i=\sum_{j\ne i}c_{ij}$, so that $w_i$ is half the degree of $i$ in the multigraph with multiplicities $2c_{ij}$. On the torus $|u_1|=\dots=|u_r|=1$ we have $u_j^{-1}=\bar u_j$, hence $(u_i-u_j)^2/(u_iu_j)=(u_i-u_j)(\bar u_j-\bar u_i)=-|u_i-u_j|^2$. Since $\prod_{i<j}(u_iu_j)^{c_{ij}}=u^w$, this gives
\[
u^{-w}\prod_{i<j}(u_i-u_j)^{2c_{ij}}=(-1)^{\sum c_{ij}}\prod_{i<j}|u_i-u_j|^{2c_{ij}} .
\]
For a Laurent polynomial $Q$, the constant term equals the average over the torus,
\[
\CT Q=\frac1{(2\pi)^r}\int_0^{2\pi}\!\!\cdots\int_0^{2\pi}Q\bigl(e^{i\theta_1},\dots,e^{i\theta_r}\bigr)\,d\theta_1\cdots d\theta_r ,
\]
because $\int_0^{2\pi}e^{ia\theta}\,d\theta=0$ for every integer $a\ne0$, so every nonconstant monomial integrates to zero. Apply this to $Q=u^{-w}\prod_{i<j}(u_i-u_j)^{2c_{ij}}$. By the previous display, $(-1)^{\sum c_{ij}}Q$ is, on the torus, a continuous nonnegative function that is positive on the nonempty open set where the $u_i$ are pairwise distinct. Hence its average is positive, and $\CT Q\ne0$. Thus the coefficient of $u^w$ in $\prod_{i<j}(u_i-u_j)^{2c_{ij}}$ is nonzero, and the multigraph with multiplicities $2c_{ij}$ has Alon--Tarsi number at most $\max_iw_i+1$. Removing edges does not increase the Alon--Tarsi number: if all monomials of $F$ with exponents at most $d$ vanish, the same holds for $(u_i-u_j)F$. Since $\mu_{ij}\le2c_{ij}$, the lemma follows.
\end{proof}

The following is the first part of Theorem~\ref{thm:intro-edges}; the rest is Corollaries~\ref{cor:extremal} and \ref{cor:random} below.

\begin{theorem}\label{thm:edges}
Let $G$ have a proper $\chi$-colouring with classes $V_1,\dots,V_\chi$, where $\chi=\chi(G)\ge2$. Then
\[
f(G)\ \le\ \max_{1\le i\le\chi}\Bigl\lfloor\frac{e(V_i,V(G)\setminus V_i)}{2}\Bigr\rfloor .
\]
In particular $f(G)\le\lfloor e(G)/2\rfloor$ and $f(G)\le\lfloor\Delta(G)\alpha(G)/2\rfloor$.
\end{theorem}

When $\chi=1$, the graph is edgeless and $f(G)=0$, so the same bound holds trivially.

\begin{proof}
Let $G'$ be the graph of Corollary~\ref{cor:embed}(b), so that $1\le\mu'_{ij}\le e_{ij}+1$, where $e_{ij}=e(V_i,V_j)$. By Corollary~\ref{cor:embed}(b) and Lemma~\ref{lem:halfdeg},
\[
f(G)\le\AT(M_{G'})-\chi\le\max_i\sum_{j\ne i}\Bigl\lceil\frac{e_{ij}+1}{2}\Bigr\rceil+1-\chi\le\max_i\sum_{j\ne i}\frac{e_{ij}}{2},
\]
using $\lceil(e+1)/2\rceil\le(e+2)/2$. The last statement holds because $e(V_i,V(G)\setminus V_i)$ is at most $e(G)$ and at most $|V_i|\Delta(G)$, and $|V_i|\le\alpha(G)$.
\end{proof}

For complete multipartite graphs with equal parts, Theorem~\ref{thm:edges} gives $f(K_{m\star\chi})\le(\chi-1)m^2/2$, which agrees with Theorem~\ref{thm:equal} to leading order.

\subsection*{The extremal function}
Following the DP literature \cite{BKZ,ZD}, let $F(n)=\max\{f(G):|V(G)|=n\}$.

\begin{corollary}\label{cor:extremal}
For every $n\ge2$,
\[
\Bigl\lfloor\frac{(\lceil n/2\rceil-1)(\lfloor n/2\rfloor-1)}{2}\Bigr\rfloor\ \le\ F(n)\ \le\ \Bigl\lfloor\frac{n^2}{8}\Bigr\rfloor .
\]
In particular $F(n)=n^2/8+O(n)$.
\end{corollary}

\begin{proof}
The lower bound is $f(K_{\lceil n/2\rceil,\lfloor n/2\rfloor})$, computed in Theorem~\ref{thm:bipartite}. For the upper bound, $e(V_i,V(G)\setminus V_i)\le|V_i|(n-|V_i|)\le n^2/4$; apply Theorem~\ref{thm:edges}.
\end{proof}

In comparison, $Z_{DP}(n)\ge n^2/4-O(n)$ \cite{BKZ}. In the worst case, therefore, the Alon--Tarsi number requires at most about half as many universal vertices as the DP-chromatic number, and the balanced complete bipartite graph is asymptotically extremal.

\subsection*{Random graphs}

\begin{corollary}\label{cor:random}
With high probability, $f(G(n,1/2))=(\frac12+o(1))\,n\log_2n$.
\end{corollary}

\begin{proof}
With high probability $e(G)=(1+o(1))n^2/4$, $\Delta(G)=(1+o(1))n/2$, $\alpha(G)=(2+o(1))\log_2n$ and $\chi(G)=(1+o(1))n/(2\log_2n)$; see \cite{BollobasBook} for the first three and \cite{Bollobas} for the last. The lower bound of Theorem~\ref{thm:main} with $G'=G$ gives $f\ge e(G)/\chi-n\ge(\frac12-o(1))n\log_2n$, and Theorem~\ref{thm:edges} gives $f\le\Delta\alpha/2\le(\frac12+o(1))n\log_2n$.
\end{proof}

\subsection*{Graphs with few vertices}
Ohba's conjecture, proved by Noel, Reed and Wu \cite{NRW}, states that every graph $G$ with $|V(G)|\le2\chi(G)+1$ is chromatic-choosable. For the paint number this fails for some graphs with $|V(G)|=2\chi(G)+1$ \cite{KKLZ}, and the on-line version of Ohba's conjecture \cite{HWZ,KKLZ} asserts that $|V(G)|\le2\chi(G)$ implies $\chi_P(G)=\chi(G)$. For general graphs, Kozik, Micek and Zhu \cite{KMZ} proved this under the stronger assumption $|V(G)|\le\chi(G)+\sqrt{\chi(G)}$. We now determine, up to an additive constant, how many vertices a graph with $\AT(G)>\chi(G)$ must have. The first step is a bound for complete multipartite graphs in terms of the number of vertices outside one vertex per class.

\begin{lemma}\label{lem:fewparts}
Let $K=K_{m_1,\dots,m_q}$ with $m_i\ge1$, and put $a_i=m_i-1$ and $t=a_1+\dots+a_q$. Then
\[
f(K)\le\max_i\Bigl\lfloor\frac{a_i(t-a_i)}2\Bigr\rfloor .
\]
\end{lemma}

\begin{proof}
For $q=1$ the graph is edgeless and both sides are $0$. For $q\ge2$, the graph $K$ is uniquely $q$-colourable, and $M_K$ has multiplicities $\mu_{ij}=m_im_j-m_i-m_j+2=a_ia_j+1$. By Theorem~\ref{thm:collapse} and Lemma~\ref{lem:halfdeg}, using $\lceil(a_ia_j+1)/2\rceil\le a_ia_j/2+1$,
\[
f(K)=\AT(M_K)-q\le\max_i\sum_{j\ne i}\Bigl(\frac{a_ia_j}2+1\Bigr)+1-q=\max_i\frac{a_i(t-a_i)}2 .
\]
Since $f(K)$ is an integer, the lemma follows.
\end{proof}

For $K_{m,m}$ the lemma gives $\lfloor(m-1)^2/2\rfloor$, which is the exact value by Theorem~\ref{thm:equal}.

\begin{theorem}\label{thm:few}
Let $G$ be a graph with $\chi(G)=\chi$ and
\[
|V(G)|\le\chi+2\sqrt{2\chi-2}-2 .
\]
Then $\AT(G)=\chi$. In particular $\chi_P(G)=\ch(G)=\chi$.
\end{theorem}

\begin{proof}
Put $t=|V(G)|-\chi$, and fix a proper $\chi$-colouring of $G$. Let $V_1,\dots,V_q$ be its classes with at least two vertices, let $s=\chi-q$ be the number of singleton classes, and put $a_i=|V_i|-1\ge1$, so that $a_1+\dots+a_q=t$. If $q=0$, then $G=K_\chi$ and $\AT(G)=\chi$. Otherwise let $H=G[V_1\cup\dots\cup V_q]$. Then $G$ is a spanning subgraph of $H\vee K_s$, where $K_s$ is formed by the singleton classes, and $\chi(H)=q$, because $\chi\le\chi(H)+s$. Also $H\subseteq K=K_{|V_1|,\dots,|V_q|}$ and $\chi(K)=q$. Adding edges does not decrease the Alon--Tarsi number (see the end of the proof of Lemma~\ref{lem:halfdeg}). Hence, by Lemma~\ref{lem:gap}, $\AT(G)\le\AT(H\vee K_s)=q+s=\chi$ as soon as $s\ge f(K)$. Since $\AT(G)\ge\chi$, it suffices to prove $s\ge f(K)$.

Let $a=\max_ia_i$. The other $q-1$ values $a_j$ are positive and sum to $t-a$, so $q\le1+t-a$ and $s\ge\chi-1-t+a$. Moreover every $a_j$ with $j$ different from an index of the maximum satisfies $a_j\le\min\{a,t-a\}$, and $x\mapsto x(t-x)$ is increasing on $[0,t/2]$ and symmetric about $t/2$, so $a_j(t-a_j)\le a(t-a)$ for all $j$. By Lemma~\ref{lem:fewparts}, $f(K)\le a(t-a)/2$. It therefore suffices that
\[
2(\chi-1-t)\ge a(t-2-a).
\]
By the inequality $xy\le\frac14(x+y)^2$, the right-hand side is at most $(t-2)^2/4$, and $2(\chi-1-t)\ge(t-2)^2/4$ is equivalent to $(t+2)^2\le8(\chi-1)$, that is, to the hypothesis.
\end{proof}

The bound of Theorem~\ref{thm:few} cannot be improved beyond an additive constant.

\begin{proposition}\label{prop:fewsharp}
Let $x\ge2$ and $2\le\chi\le\lfloor x^2/2\rfloor+1$. Then $G=K_{x+1,x+1}\vee K_{\chi-2}$ satisfies $\chi(G)=\chi$, $|V(G)|=\chi+2x$ and $\AT(G)>\chi$.
\end{proposition}

\begin{proof}
By Theorem~\ref{thm:equal}, $f(K_{x+1,x+1})=\lfloor x^2/2\rfloor>\chi-2$. Since $f(K_{x+1,x+1})$ is the least $k$ with $\AT(K_{x+1,x+1}\vee K_k)=2+k$, this gives $\AT(G)\ne2+(\chi-2)=\chi$. As $\AT(G)\ge\chi(G)=\chi$, we conclude that $\AT(G)>\chi$.
\end{proof}

Let $\psi_{\AT}(\chi)$ be the least number of vertices of a graph $G$ with $\chi(G)=\chi$ and $\AT(G)>\chi$. Taking $x=\lceil\sqrt{2\chi-1}\rceil$ in Proposition~\ref{prop:fewsharp}, so that $\lfloor x^2/2\rfloor\ge\chi-1$, we obtain:

\begin{corollary}\label{cor:psi}
For every $\chi\ge2$,
\[
\chi+2\sqrt{2\chi-2}-2<\psi_{\AT}(\chi)<\chi+2\sqrt{2\chi-1}+2 .
\]
In particular $\psi_{\AT}(\chi)=\chi+2\sqrt{2\chi}+O(1)$.
\end{corollary}

Since $\chi_P\le\AT$, Theorem~\ref{thm:few} improves the bound $|V(G)|\le\chi(G)+\sqrt{\chi(G)}$ of Kozik, Micek and Zhu for the paint number for every $\chi\ge3$. Corollary~\ref{cor:psi} also shows the limits of this approach. The on-line version of Ohba's conjecture concerns graphs with up to $2\chi$ vertices, but already for $|V(G)|=\chi+(2\sqrt2+o(1))\sqrt\chi$ there are graphs with $\AT(G)>\chi$. So the Alon--Tarsi number alone cannot prove that conjecture, nor any bound of the form $|V(G)|\le\chi+C\sqrt\chi$ with $C>2\sqrt2$.

The same comparison applies to joins. By the definition of $f$, the graph $G\vee K_k$ is chromatic-paintable for every $k\ge f(G)$, so, by Corollary~\ref{cor:extremal}, $k\ge\lfloor n^2/8\rfloor$ always suffices, and Corollary~\ref{cor:random} gives $(\frac12+o(1))n\log_2n$ for $G(n,1/2)$. The general bound of \cite{KMZ}, applied to $G\vee K_k$, requires $k\ge(n-\chi)^2-\chi$.

\section{Bipartite graphs}\label{sec:bip}

A bipartite graph need not be uniquely colourable, but each connected component has exactly two proper maps into a fixed two-colour set. A component containing an edge uses both colours in each map; an isolated vertex does not. The colouring polynomial nevertheless factorises over the components.

\begin{proof}[Proof of Theorem~\ref{thm:bipartite}]
Here $\chi=2$ and $D=k+1$. Let $C$ range over the components of $G$, with sides $X_C,Y_C$, $n_C$ vertices and $e_C$ edges, and call $C$ \emph{odd} if $e_C-n_C$ is odd. Write $S=\{s,t\}$ for the two colours. Each component $C$ has exactly two maps into $S$ that are proper on $C$: one sends $X_C$ to $s$ and $Y_C$ to $t$, the other does the reverse. The first has deficit weight $A_C(s,t)=s^{\sum_{v\in X_C}r_v}t^{\sum_{v\in Y_C}r_v}$, and the second has weight $A_C(t,s)$. For an isolated vertex, take $X_C=\{v\}$ and $Y_C=\varnothing$; its two maps into $S$ give, up to the orientation sign $\eta_C$, the factor $(s^{r_v}-t^{r_v})/(s-t)$. Thus the factorisation $\Phi_{G,r}=\prod_CF_C$ is valid for all components, where
\[
F_C=\eta_C\Bigl[(s-t)^{e_C-n_C}A_C(s,t)+(t-s)^{e_C-n_C}A_C(t,s)\Bigr],
\]
where the sign $\eta_C=\pm1$ depends only on the vertex order within $C$. The exponent $e_C-n_C$ can equal $-1$ (for a tree); after combining the two terms, the numerator $A_C(s,t)-A_C(t,s)$ is divisible by $s-t$, so $F_C$ is still a polynomial. For an isolated vertex the factor reduces to the divided difference just described.

\emph{Lower bound.} If $C$ is odd and $r$ vanishes on $C$, then $F_C=0$. By Proposition~\ref{prop:reduction}, a nonzero admissible monomial must therefore have a positive deficit in each odd component, and \eqref{eq:budget} gives $n+2k-e(G)\ge r_{\rm odd}$, the number of odd components. Hence
\[
k\ge\frac{e(G)-n+r_{\rm odd}}{2}=\sum_C\Bigl\lceil\frac{e_C-n_C}{2}\Bigr\rceil=\sum_C\Bigl\lfloor\frac{\beta(C)}2\Bigr\rfloor.
\]

\emph{Upper bound.} Let $k$ be equal to this sum. Then the budget equals $r_{\rm odd}$. For each odd component containing an edge, give one vertex of $Y_C$ deficit $1$; for an isolated component, give its single vertex deficit $1$. Give all other vertices deficit $0$, and take $\lambda=(2^k)$. Then $F_C=\pm2(s-t)^{e_C-n_C}$ for even $C$, and $F_C=\pm(s-t)^{e_C-n_C}(t-s)$ for odd $C$. Hence $\Phi_{G,r}(s,t)=c(s-t)^{2k}$ with $c\ne0$, and by Proposition~\ref{prop:reduction} and Lemma~\ref{lem:schur}, applied with $\lambda'=(k,k)$, the coefficient of the corresponding monomial is
\[
\pm c\,\coef{s^{k+1}t^{k}}(s-t)^{2k+1}=\pm c\binom{2k+1}{k}\ne0.\qedhere
\]
\end{proof}

For connected $G$, Theorem~\ref{thm:bipartite} also follows from Theorem~\ref{thm:collapse}: a connected bipartite graph is uniquely $2$-colourable, $M_G$ consists of $\beta(G)+1$ parallel edges, and a two-vertex multigraph with $\mu$ edges has Alon--Tarsi number $\lceil\mu/2\rceil+1$, because the coefficients of $(u-v)^\mu$ are $\pm\binom\mu j$, all nonzero, and some exponent of each monomial is at least $\mu/2$. For $K_{p,q}$ we obtain $f(K_{p,q})=\lfloor(p-1)(q-1)/2\rfloor$.

\section{Complete multipartite graphs: a conjecture}\label{sec:multi}

For $\vec m=(m_1,\dots,m_\chi)$ let $M^+(\vec m)=M_{K_{\vec m}}$; thus $\mu_{ij}=(m_i-1)(m_j-1)+1$. By Theorem~\ref{thm:collapse}, $f(K_{\vec m})=\AT(M^+(\vec m))-\chi$, and every lower bound for the Alon--Tarsi number of $M^+(\vec m)$ yields one for $f$. The bound below combines two obstructions, edge density and twins, and improves on Proposition~\ref{prop:density}.

Throughout this section and the next, $M$ is a loopless multigraph on $\{1,\dots,\chi\}$ with positive multiplicity $\mu_{ij}$ between any two distinct vertices; $M^+(\vec m)$ is of this kind. Since $M$ contains $K_\chi$, we have $\AT(M)\ge\chi$.

\begin{definition}
Two vertices $i\ne j$ of $M$ are \emph{twins} if $\mu_{ij}$ is odd and $\mu_{il}=\mu_{jl}$ for every $l\notin\{i,j\}$. (As for the twins of Lemma~\ref{lem:twin}, interchanging $u_i$ and $u_j$ then changes the sign of $P_M$.)
\end{definition}

If $i,j$ and $j,l$ are twins, then $\mu_{ij}=\mu_{il}=\mu_{jl}$ and $i,l$ are twins as well. Hence the relation ``$i=j$ or $i,j$ are twins'' is an equivalence relation; we call its classes the \emph{twin classes} of $M$. For $M=M^+(\vec m)$ we have $\mu_{il}-\mu_{jl}=(m_i-m_j)(m_l-1)$, so the condition $\mu_{il}=\mu_{jl}$ for all $l\notin\{i,j\}$ holds exactly when $m_i=m_j$ or all other parts have size $1$. For instance, two parts of the same odd size are twins, and for $\vec m=(1,2,3)$ the parts of sizes $2$ and $3$ are twins.

\begin{proposition}\label{prop:obstructions}
Let $J\subseteq\{1,\dots,\chi\}$ be nonempty, put $e_J=\sum_{\{i,j\}\subseteq J}\mu_{ij}$, and let $t_1,\dots,t_s$ be the sizes of the nonempty intersections of $J$ with the twin classes of $M$. If $P_M$ has a nonzero monomial with all exponents at most $D$, then
\begin{equation}\label{eq:twindensity}
\sum_{l=1}^s\Bigl(t_lD-\binom{t_l}{2}\Bigr)\ \ge\ e_J .
\end{equation}
\end{proposition}

\begin{proof}
In every monomial of $P_M$ the exponents of the variables indexed by $J$ sum to at least $e_J$, since each edge inside $J$ contributes one of its endpoints. If $i,j$ are twins, then interchanging $u_i$ and $u_j$ multiplies $P_M$ by $(-1)^{\mu_{ij}}=-1$, so as in Lemma~\ref{lem:twin} they must receive distinct exponents. The $t$ variables of $J$ in one twin class therefore have distinct exponents in $[0,D]$, whose sum is at most $tD-\binom t2$.
\end{proof}

\begin{definition}
Let $D_0(M)$ be the least $D\ge\chi-1$ for which \eqref{eq:twindensity} holds for every nonempty $J$, and write $D_0(\vec m)=D_0(M^+(\vec m))$.
\end{definition}

Since $t_1+\dots+t_s=|J|$, inequality \eqref{eq:twindensity} can be rewritten as $|J|\,D\ge e_J+\sum_l\binom{t_l}2$. Hence
\begin{equation}\label{eq:D0formula}
D_0(M)=\max\biggl\{\chi-1,\ \max_{\varnothing\ne J}\biggl\lceil\frac{e_J+\sum_{C}\binom{|C\cap J|}{2}}{|J|}\biggr\rceil\biggr\},
\end{equation}
where $C$ runs over the twin classes of $M$. The first term in the fraction gives the induced edge count per vertex in $J$, and the second adds $1/|J|$ for every pair of twins in $J$. The floor $\chi-1$ comes from $\AT(M)\ge\chi$; it also ensures that each twin class, which has at most $\chi$ vertices, can receive distinct exponents in $[0,D]$.

By Proposition~\ref{prop:obstructions} and since $\AT(M)\ge\chi$, we get $\AT(M)\ge D_0(M)+1$. We call this, and the resulting bound on $f$, the \emph{twin-density bound}. Together with Theorem~\ref{thm:collapse} it gives
\[
f(K_{\vec m})\ \ge\ D_0(\vec m)+1-\chi .
\]

\begin{example}
For $K_{3,3,3}$ the density bound of Proposition~\ref{prop:density} gives $3$, whereas $D_0=6$ and so $f\ge4$; by Theorem~\ref{thm:equal}, $f(K_{3,3,3})=4$. Here the gain comes from twins: all three classes are mutual twins. For $K_{2,5,8}$ there are no twins, and the gain comes from measuring density on a subset of the classes rather than on the whole graph: Proposition~\ref{prop:density} gives $12$, whereas $J=\{2,3\}$ gives $D_0\ge\lceil29/2\rceil=15$ and so $f\ge13$, in agreement with Example~\ref{ex:gaps}. In fact $f(K_{2,5,8})=13$: here $\mu=(5,8,29)$, $e=42$ and $D^*=15$, so $3D^*>e$ and Theorem~\ref{thm:three} below gives $\AT(M^+)=16$.
\end{example}

\begin{conjecture}\label{conj:M}
For every $\vec m$, $\AT(M^+(\vec m))=D_0(\vec m)+1$. Equivalently (recall that $D_0(\vec m)\ge\chi-1$),
\[
f(K_{\vec m})=D_0(\vec m)+1-\chi.
\]
\end{conjecture}

The conjecture says that, for complete multipartite graphs, edge density and twins are the only obstructions to a nonzero coefficient with small exponents. This is a statement about the special multiplicities $(m_i-1)(m_j-1)+1$: for general multigraphs other cancellations occur, and Example~\ref{ex:3jzero} gives a three-vertex multigraph with $\AT(M)=D_0(M)+2$.

Theorems~\ref{thm:equal} and \ref{thm:bipartite} prove Conjecture~\ref{conj:M} when all parts have equal size and when $\chi=2$; Section~\ref{sec:three} extends this considerably. We have verified it by computer for all $\vec m$ with $\chi=3$ and $m_i\le8$, with $\chi=4$ and $m_i\le5$, and with $\chi=5$ and $m_i\le3$ (211 vectors). If the conjecture holds, then every graph $G$ with a proper $\chi$-colouring with class sizes $\vec m$ satisfies $f(G)\le D_0(\vec m)+1-\chi$, a bound that is exact for $K_{\vec m}$.

Table~\ref{tab:data} lists values of $f$ obtained by a direct computation that does not use Sections~\ref{sec:reduction}--\ref{sec:bip}. For $k=0,1,2,\dots$ we enumerated the admissible exponent vectors and evaluated their coefficients exactly, in rational arithmetic, with Theorem~\ref{thm:qcn} and the grids $\{0,1,\dots,a_v\}$, stopping at the first $k$ with a nonzero coefficient. Only terms in which $G$ is properly coloured and the universal vertices receive distinct colours not used on $G$ are nonzero, which keeps the sums small, and below the density bound of Proposition~\ref{prop:density} there are no admissible vectors at all. All values agree with Theorems~\ref{thm:collapse}, \ref{thm:bipartite} and \ref{thm:equal}.

\begin{table}[ht]
\centering\small
\begin{tabular}{lcc@{\qquad}lcc@{\qquad}lcc}
\toprule
$G$ & (P) & $f$ & $G$ & (P) & $f$ & $G$ & (P) & $f$\\
\midrule
$K_{3,3}$&2&2 & $K_{2,2,3}$&1&1 & $2K_{2,3}$&1&2\\
$K_{4,4}$&4&4 & $K_{2,3,3}$&2&2 & $K_{3,3}\sqcup K_{2,3}$&2&3\\
$K_{5,5}$&8&8 & $K_{2,2,4}$&2&2 & $K_{3,3}-e$&1&1\\
$K_{6,6}$&12&12 & $K_{2,3,4}$&3&3 & $Q_3$&2&2\\
$K_{2,5}$&2&2 & $K_{2,4,4}$&4&4 & $K_{2,3,3}-e$&2&2\\
$K_{3,5}$&4&4 & $K_{3,3,4}$&5&5 & $K_{3,3,3}-e$&3&3\\
$K_{4,6}$&7&7 & $K_{3,4,4}$&6&6 & $K_{4,4,4}$&8&8\\
$K_{2,2,2,2}$&0&0 & $K_{3,3,3}$&3&4 & $K_{3,3,5}$&6&6\\
\bottomrule
\end{tabular}
\caption{Values of $f$ computed directly from Theorem~\ref{thm:qcn}, compared with the lower bound (P) of Proposition~\ref{prop:density} applied to $G'=G$.}\label{tab:data}
\end{table}

\section{The conjecture for three classes and for one exceptional class}\label{sec:three}

We prove Conjecture~\ref{conj:M} in two cases. For $\chi=3$ we determine $\AT(M)$ for every three-vertex multigraph $M$ by a three-term recurrence for the coefficients of $P_M$ (Theorem~\ref{thm:three}). For $\vec m=(p,q,\dots,q)$ the relevant coefficient is given by Morris's constant term identity (Theorem~\ref{thm:morris}).

\subsection{Three-vertex multigraphs}
Let $M$ be a multigraph on $\{1,2,3\}$ with positive multiplicities $\mu_{12},\mu_{13},\mu_{23}$, and write
\[
P=P_M=(u_1-u_2)^{\mu_{12}}(u_1-u_3)^{\mu_{13}}(u_2-u_3)^{\mu_{23}}=\sum_\gamma p_\gamma u^\gamma,\qquad e=\mu_{12}+\mu_{13}+\mu_{23},
\]
with vertex degrees $d_1=\mu_{12}+\mu_{13}$, $d_2=\mu_{12}+\mu_{23}$, $d_3=\mu_{13}+\mu_{23}$; we put $p_\gamma=0$ when $u^\gamma$ is not a monomial of $P$. Every monomial lies in $R=\{\gamma\in\mathbb Z^3:\sum\gamma_k=e,\ 0\le\gamma_k\le d_k\}$. Let $\varepsilon_1,\varepsilon_2,\varepsilon_3$ be the standard basis vectors of $\mathbb Z^3$; we call $\gamma,\gamma'$ \emph{adjacent} if $\gamma-\gamma'=\varepsilon_i-\varepsilon_j$ for some $i\ne j$.

The proof uses three first-order differential operators that annihilate $P$. Let $\partial_k=\partial/\partial u_k$, and for $k=1,2,3$ put
\[
E_k=\partial_k,\qquad F_k=u_k^2\partial_k-d_ku_k,\qquad H_k=u_k\partial_k-\tfrac12d_k,
\]
and $E=\sum_kE_k$, $F=\sum_kF_k$, $H=\sum_kH_k$. They come from three symmetries of $P$: $E$ from invariance under the translation $u\mapsto u+t\one$, $H$ from homogeneity, and $F$ from the identity $P(u_1^{-1},u_2^{-1},u_3^{-1})\,u_1^{d_1}u_2^{d_2}u_3^{d_3}=\pm P(u)$, which holds because $u_i^{-1}-u_j^{-1}=(u_j-u_i)/(u_iu_j)$.

\begin{lemma}\label{lem:sl2}
$EP=FP=HP=0$. Equivalently, $\sum_k(\tau_k+1)p_{\tau+\varepsilon_k}=0$ whenever $\sum\tau_k=e-1$, and $\sum_k(\tau_k-1-d_k)p_{\tau-\varepsilon_k}=0$ whenever $\sum\tau_k=e+1$.
\end{lemma}

\begin{proof}
$EP=0$ because $P(u+t\one)=P(u)$, and $HP=0$ because $P$ is homogeneous of degree $e=\frac12\sum d_k$. For $F$, note that $(u_i^2\partial_i+u_j^2\partial_j)(u_i-u_j)=(u_i+u_j)(u_i-u_j)$, so $\sum_ku_k^2\partial_kP=P\sum_{i<j}\mu_{ij}(u_i+u_j)=P\sum_kd_ku_k$. The two identities are the coefficients of $u^\tau$ in $EP=0$ and $FP=0$.
\end{proof}

\begin{lemma}[three-term recurrence]\label{lem:rec}
For every $\gamma\in\mathbb Z^3$ with $\sum\gamma_k=e$,
\begin{multline*}
2(\gamma_1+1)(d_2+1-\gamma_2)\,p_{\gamma+\varepsilon_1-\varepsilon_2}+2(d_1+1-\gamma_1)(\gamma_2+1)\,p_{\gamma-\varepsilon_1+\varepsilon_2}\\
=\bigl[\mu_{13}\mu_{23}-\mu_{12}(e+2)-(2\gamma_1-d_1)(2\gamma_2-d_2)\bigr]p_\gamma .
\end{multline*}
\end{lemma}

\begin{proof}
For each $k$ let $\Omega_k=2H_k^2-E_kF_k-F_kE_k$. On a monomial $u_k^j$,
\[
E_kF_ku_k^j=(j-d_k)(j+1)\,u_k^j,\qquad F_kE_ku_k^j=j(j-1-d_k)\,u_k^j,\qquad 2H_k^2u_k^j=2\bigl(j-\tfrac12d_k\bigr)^2u_k^j,
\]
so $\Omega_k$ acts on all polynomials in $u_k$ as multiplication by $\theta_k=d_k+\frac12d_k^2$. Let $\Omega_{12}$ be the same expression with $E_k,F_k,H_k$ replaced by $E_1+E_2$, $F_1+F_2$, $H_1+H_2$. By Lemma~\ref{lem:sl2}, $(E_1+E_2)P=-E_3P$, and likewise for $F$ and $H$. Since operators in different variables commute, it follows that $(H_1+H_2)^2P=H_3^2P$, $(E_1+E_2)(F_1+F_2)P=F_3E_3P$ and $(F_1+F_2)(E_1+E_2)P=E_3F_3P$, so $\Omega_{12}P=\Omega_3P=\theta_3P$. On the other hand, expanding gives $\Omega_{12}=\Omega_1+\Omega_2+4H_1H_2-2E_1F_2-2E_2F_1$, so
\[
(4H_1H_2-2E_1F_2-2E_2F_1)P=(\theta_3-\theta_1-\theta_2)P,\qquad \theta_3-\theta_1-\theta_2=\mu_{13}\mu_{23}-\mu_{12}(e+2).
\]
Finally compare coefficients of $u^\gamma$: $H_1H_2$ multiplies $u^\gamma$ by $(\gamma_1-\frac12d_1)(\gamma_2-\frac12d_2)$, while $E_1F_2$ maps $u^{\gamma+\varepsilon_1-\varepsilon_2}$ to $(\gamma_1+1)(\gamma_2-1-d_2)u^\gamma$ and $E_2F_1$ maps $u^{\gamma-\varepsilon_1+\varepsilon_2}$ to $(\gamma_2+1)(\gamma_1-1-d_1)u^\gamma$.
\end{proof}

For $\gamma\in R$ both coefficients on the left-hand side of Lemma~\ref{lem:rec} are nonzero.

\begin{proposition}\label{prop:noadj}
\leavevmode
\begin{enumerate}[label=\textup{(\alph*)}]
\item For every $0\le t\le d_3$ some $\gamma\in R$ with $\gamma_3=t$ has $p_\gamma\ne0$.
\item If $\gamma,\gamma'\in R$ are adjacent, then $p_\gamma$ and $p_{\gamma'}$ are not both zero.
\end{enumerate}
\end{proposition}

\begin{proof}
(a) Suppose that $p_\gamma=0$ whenever $\gamma_3=t$. For $\gamma\in R$ with $\gamma_3=t+1$, the first identity of Lemma~\ref{lem:sl2} with $\tau=\gamma-\varepsilon_3$ expresses $(t+1)p_\gamma$ through coefficients with third exponent $t$, so $p_\gamma=0$; inductively all coefficients with third exponent above $t$ vanish. The second identity with $\tau=\gamma+\varepsilon_3$ handles third exponent $t-1$ in the same way, with the nonzero weight $t-1-d_3$, and inductively everything below $t$. Then $P=0$, which is absurd.
(b) By symmetry let $\gamma'=\gamma+\varepsilon_1-\varepsilon_2$ and $t=\gamma_3$. On the slice $L=\{x\in R:x_3=t\}$ both outer coefficients of Lemma~\ref{lem:rec} are nonzero, so two consecutive zeros propagate along $L$ in both directions, and $p_x=0$ for all $x\in L$, contradicting (a).
\end{proof}

\begin{figure}[ht]
\centering
\begin{tikzpicture}[scale=0.42]
\fill[blue!12] (13.500,7.794) -- (12.000,5.196) -- (10.500,7.794) -- cycle;
\fill[black!70] (1.000,1.732) circle (0.14);
\fill[black!70] (1.500,2.598) circle (0.14);
\fill[black!70] (2.000,3.464) circle (0.14);
\fill[black!70] (2.500,4.330) circle (0.14);
\fill[black!70] (3.000,5.196) circle (0.14);
\fill[black!70] (3.500,6.062) circle (0.14);
\fill[black!70] (4.000,6.928) circle (0.14);
\fill[black!70] (4.500,7.794) circle (0.14);
\fill[black!70] (5.000,8.660) circle (0.14);
\fill[black!70] (5.500,9.526) circle (0.14);
\fill[black!70] (6.000,10.392) circle (0.14);
\fill[black!70] (6.500,11.258) circle (0.14);
\fill[black!70] (7.000,12.124) circle (0.14);
\fill[black!70] (7.500,12.990) circle (0.14);
\fill[black!70] (1.500,0.866) circle (0.14);
\fill[black!70] (2.000,1.732) circle (0.14);
\fill[black!70] (2.500,2.598) circle (0.14);
\fill[black!70] (3.000,3.464) circle (0.14);
\fill[black!70] (3.500,4.330) circle (0.14);
\fill[black!70] (4.000,5.196) circle (0.14);
\fill[black!70] (4.500,6.062) circle (0.14);
\fill[black!70] (5.000,6.928) circle (0.14);
\fill[black!70] (5.500,7.794) circle (0.14);
\fill[black!70] (6.000,8.660) circle (0.14);
\fill[black!70] (6.500,9.526) circle (0.14);
\fill[black!70] (7.000,10.392) circle (0.14);
\fill[black!70] (7.500,11.258) circle (0.14);
\fill[black!70] (8.000,12.124) circle (0.14);
\fill[black!70] (8.500,12.990) circle (0.14);
\fill[black!70] (2.000,0.000) circle (0.14);
\fill[black!70] (2.500,0.866) circle (0.14);
\fill[black!70] (3.000,1.732) circle (0.14);
\fill[black!70] (3.500,2.598) circle (0.14);
\fill[black!70] (4.000,3.464) circle (0.14);
\fill[black!70] (4.500,4.330) circle (0.14);
\fill[black!70] (5.000,5.196) circle (0.14);
\fill[black!70] (5.500,6.062) circle (0.14);
\fill[black!70] (6.000,6.928) circle (0.14);
\fill[black!70] (6.500,7.794) circle (0.14);
\fill[black!70] (7.000,8.660) circle (0.14);
\fill[black!70] (7.500,9.526) circle (0.14);
\fill[black!70] (8.000,10.392) circle (0.14);
\fill[black!70] (8.500,11.258) circle (0.14);
\fill[black!70] (9.000,12.124) circle (0.14);
\fill[black!70] (9.500,12.990) circle (0.14);
\fill[black!70] (3.000,0.000) circle (0.14);
\fill[black!70] (3.500,0.866) circle (0.14);
\fill[black!70] (4.000,1.732) circle (0.14);
\fill[black!70] (4.500,2.598) circle (0.14);
\fill[black!70] (5.000,3.464) circle (0.14);
\fill[black!70] (5.500,4.330) circle (0.14);
\fill[black!70] (6.000,5.196) circle (0.14);
\draw[red,thick,fill=white] (6.500,6.062) circle (0.32);
\fill[black!70] (7.000,6.928) circle (0.14);
\fill[black!70] (7.500,7.794) circle (0.14);
\fill[black!70] (8.000,8.660) circle (0.14);
\fill[black!70] (8.500,9.526) circle (0.14);
\fill[black!70] (9.000,10.392) circle (0.14);
\fill[black!70] (9.500,11.258) circle (0.14);
\fill[black!70] (10.000,12.124) circle (0.14);
\fill[black!70] (10.500,12.990) circle (0.14);
\fill[black!70] (4.000,0.000) circle (0.14);
\fill[black!70] (4.500,0.866) circle (0.14);
\fill[black!70] (5.000,1.732) circle (0.14);
\fill[black!70] (5.500,2.598) circle (0.14);
\fill[black!70] (6.000,3.464) circle (0.14);
\fill[black!70] (6.500,4.330) circle (0.14);
\fill[black!70] (7.000,5.196) circle (0.14);
\fill[black!70] (7.500,6.062) circle (0.14);
\fill[black!70] (8.000,6.928) circle (0.14);
\fill[black!70] (8.500,7.794) circle (0.14);
\fill[black!70] (9.000,8.660) circle (0.14);
\fill[black!70] (9.500,9.526) circle (0.14);
\fill[black!70] (10.000,10.392) circle (0.14);
\fill[black!70] (10.500,11.258) circle (0.14);
\fill[black!70] (11.000,12.124) circle (0.14);
\fill[black!70] (11.500,12.990) circle (0.14);
\fill[black!70] (5.000,0.000) circle (0.14);
\fill[black!70] (5.500,0.866) circle (0.14);
\fill[black!70] (6.000,1.732) circle (0.14);
\fill[black!70] (6.500,2.598) circle (0.14);
\fill[black!70] (7.000,3.464) circle (0.14);
\fill[black!70] (7.500,4.330) circle (0.14);
\fill[black!70] (8.000,5.196) circle (0.14);
\fill[black!70] (8.500,6.062) circle (0.14);
\fill[black!70] (9.000,6.928) circle (0.14);
\fill[black!70] (9.500,7.794) circle (0.14);
\fill[black!70] (10.000,8.660) circle (0.14);
\fill[black!70] (10.500,9.526) circle (0.14);
\fill[black!70] (11.000,10.392) circle (0.14);
\fill[black!70] (11.500,11.258) circle (0.14);
\fill[black!70] (12.000,12.124) circle (0.14);
\fill[black!70] (12.500,12.990) circle (0.14);
\fill[black!70] (6.000,0.000) circle (0.14);
\fill[black!70] (6.500,0.866) circle (0.14);
\fill[black!70] (7.000,1.732) circle (0.14);
\fill[black!70] (7.500,2.598) circle (0.14);
\fill[black!70] (8.000,3.464) circle (0.14);
\fill[black!70] (8.500,4.330) circle (0.14);
\fill[black!70] (9.000,5.196) circle (0.14);
\fill[black!70] (9.500,6.062) circle (0.14);
\fill[black!70] (10.000,6.928) circle (0.14);
\fill[black!70] (10.500,7.794) circle (0.14);
\fill[black!70] (11.000,8.660) circle (0.14);
\fill[black!70] (11.500,9.526) circle (0.14);
\fill[black!70] (12.000,10.392) circle (0.14);
\fill[black!70] (12.500,11.258) circle (0.14);
\fill[black!70] (13.000,12.124) circle (0.14);
\fill[black!70] (13.500,12.990) circle (0.14);
\fill[black!70] (7.000,0.000) circle (0.14);
\fill[black!70] (7.500,0.866) circle (0.14);
\fill[black!70] (8.000,1.732) circle (0.14);
\fill[black!70] (8.500,2.598) circle (0.14);
\fill[black!70] (9.000,3.464) circle (0.14);
\fill[black!70] (9.500,4.330) circle (0.14);
\fill[black!70] (10.000,5.196) circle (0.14);
\fill[black!70] (10.500,6.062) circle (0.14);
\fill[black!70] (11.000,6.928) circle (0.14);
\fill[black!70] (11.500,7.794) circle (0.14);
\fill[black!70] (12.000,8.660) circle (0.14);
\fill[black!70] (12.500,9.526) circle (0.14);
\fill[black!70] (13.000,10.392) circle (0.14);
\fill[black!70] (13.500,11.258) circle (0.14);
\fill[black!70] (14.000,12.124) circle (0.14);
\fill[black!70] (14.500,12.990) circle (0.14);
\fill[black!70] (8.000,0.000) circle (0.14);
\fill[black!70] (8.500,0.866) circle (0.14);
\fill[black!70] (9.000,1.732) circle (0.14);
\fill[black!70] (9.500,2.598) circle (0.14);
\fill[black!70] (10.000,3.464) circle (0.14);
\fill[black!70] (10.500,4.330) circle (0.14);
\fill[black!70] (11.000,5.196) circle (0.14);
\fill[black!70] (11.500,6.062) circle (0.14);
\draw[red,thick,fill=white] (12.000,6.928) circle (0.32);
\fill[black!70] (12.500,7.794) circle (0.14);
\fill[black!70] (13.000,8.660) circle (0.14);
\fill[black!70] (13.500,9.526) circle (0.14);
\fill[black!70] (14.000,10.392) circle (0.14);
\fill[black!70] (14.500,11.258) circle (0.14);
\fill[black!70] (15.000,12.124) circle (0.14);
\fill[black!70] (15.500,12.990) circle (0.14);
\fill[black!70] (9.000,0.000) circle (0.14);
\fill[black!70] (9.500,0.866) circle (0.14);
\fill[black!70] (10.000,1.732) circle (0.14);
\fill[black!70] (10.500,2.598) circle (0.14);
\fill[black!70] (11.000,3.464) circle (0.14);
\fill[black!70] (11.500,4.330) circle (0.14);
\fill[black!70] (12.000,5.196) circle (0.14);
\fill[black!70] (12.500,6.062) circle (0.14);
\fill[black!70] (13.000,6.928) circle (0.14);
\fill[black!70] (13.500,7.794) circle (0.14);
\fill[black!70] (14.000,8.660) circle (0.14);
\fill[black!70] (14.500,9.526) circle (0.14);
\fill[black!70] (15.000,10.392) circle (0.14);
\fill[black!70] (15.500,11.258) circle (0.14);
\fill[black!70] (16.000,12.124) circle (0.14);
\fill[black!70] (16.500,12.990) circle (0.14);
\fill[black!70] (10.000,0.000) circle (0.14);
\fill[black!70] (10.500,0.866) circle (0.14);
\fill[black!70] (11.000,1.732) circle (0.14);
\fill[black!70] (11.500,2.598) circle (0.14);
\fill[black!70] (12.000,3.464) circle (0.14);
\fill[black!70] (12.500,4.330) circle (0.14);
\fill[black!70] (13.000,5.196) circle (0.14);
\fill[black!70] (13.500,6.062) circle (0.14);
\fill[black!70] (14.000,6.928) circle (0.14);
\fill[black!70] (14.500,7.794) circle (0.14);
\fill[black!70] (15.000,8.660) circle (0.14);
\fill[black!70] (15.500,9.526) circle (0.14);
\fill[black!70] (16.000,10.392) circle (0.14);
\fill[black!70] (16.500,11.258) circle (0.14);
\fill[black!70] (17.000,12.124) circle (0.14);
\fill[black!70] (11.000,0.000) circle (0.14);
\fill[black!70] (11.500,0.866) circle (0.14);
\fill[black!70] (12.000,1.732) circle (0.14);
\fill[black!70] (12.500,2.598) circle (0.14);
\fill[black!70] (13.000,3.464) circle (0.14);
\fill[black!70] (13.500,4.330) circle (0.14);
\fill[black!70] (14.000,5.196) circle (0.14);
\fill[black!70] (14.500,6.062) circle (0.14);
\fill[black!70] (15.000,6.928) circle (0.14);
\fill[black!70] (15.500,7.794) circle (0.14);
\fill[black!70] (16.000,8.660) circle (0.14);
\fill[black!70] (16.500,9.526) circle (0.14);
\fill[black!70] (17.000,10.392) circle (0.14);
\fill[black!70] (17.500,11.258) circle (0.14);
\draw[-{Stealth},thick] (-2,-1.2) -- (1,-1.2) node[right] {\scriptsize $\gamma_1$};
\draw[-{Stealth},thick] (-2,-1.2) -- (-0.5,1.398) node[above] {\scriptsize $\gamma_2$};
\draw[red!70!black,thin] (12.300,7.128) -- (19.300,9.928) node[right,fill=white,inner sep=1pt] {\scriptsize $(8,8,8)$};
\draw[red!70!black,thin] (6.200,5.862) -- (2.000,3.562) node[left,fill=white,inner sep=1pt] {\scriptsize $(3,7,14)$};
\end{tikzpicture}
\caption{The $186$ exponent vectors $\gamma\in R$ for $\mu=(2,9,13)$, drawn in the plane $\gamma_1+\gamma_2+\gamma_3=24$. Only two coefficients vanish (open circles), and they are not adjacent, as Proposition~\ref{prop:noadj}(b) requires. The shaded triangle is $R_9=R\cap[0,9]^3$. Its centre $(8,8,8)$ is the only point of $R_8$, and its coefficient vanishes, so $\AT(M)=10$ (Example~\ref{ex:3jzero}).}\label{fig:support}
\end{figure}
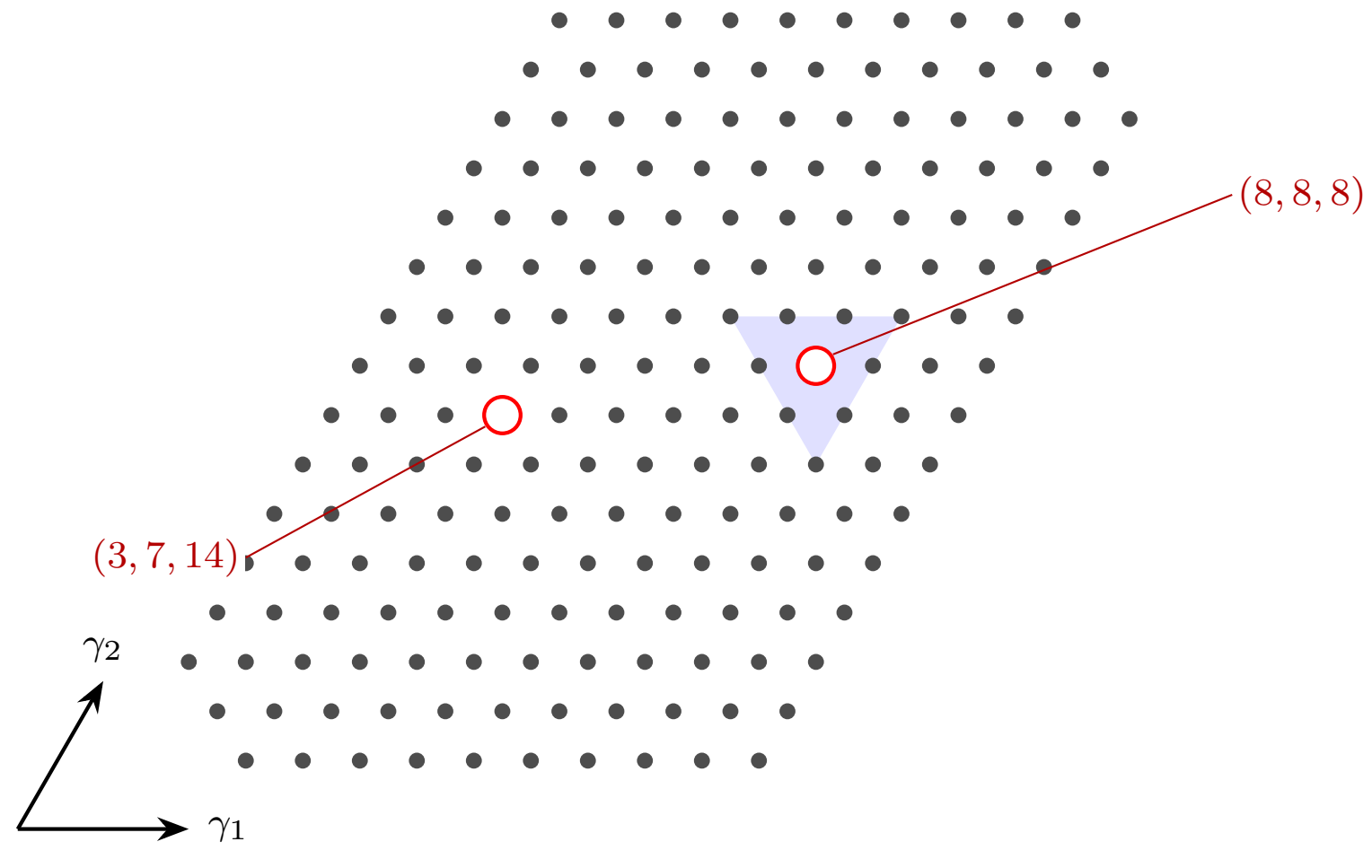

Figure~\ref{fig:support} illustrates Proposition~\ref{prop:noadj}(b) for one multigraph.

For an integer $D$ let $R_D=R\cap[0,D]^3$, and let $D^*$ be the least $D$ with $R_D\ne\emptyset$. Since $R_D\ne\emptyset$ if and only if $\sum_k\min(d_k,D)\ge e$, and $e-d_k$ is the multiplicity of the pair not containing $k$, expanding the minimum shows that $D^*$ is the least $D$ with $3D\ge e$ and $2D\ge\mu_{ij}$ for all pairs. Clearly $\AT(M)\ge D^*+1$. These are the conditions \eqref{eq:twindensity} without the twin corrections, and we shall see that $D^*=D_0(M)$ except in one case involving twins.

\begin{theorem}\label{thm:three}
Let $M$ be a multigraph on three vertices with positive multiplicities. Then $\AT(M)=D^*+1$, unless $3D^*=e$, $D^*<\min_kd_k$ and the central coefficient $p_{(D^*,D^*,D^*)}$ vanishes, in which case $\AT(M)=D^*+2$. Consequently $\AT(M)=D_0(M)+1$, unless $M$ has no twins, $3\mid e$, $D_0=e/3<\min_kd_k$ and $p_{(D_0,D_0,D_0)}=0$, in which case $\AT(M)=D_0(M)+2$.
\end{theorem}

\begin{proof}
Fix $D\ge D^*$. If $|R_D|\ge2$, take $x\ne y$ in $R_D$ and indices with $x_i>y_i$, $x_j<y_j$; then $x-\varepsilon_i+\varepsilon_j\in R_D$ is adjacent to $x$, and by Proposition~\ref{prop:noadj} one of them has a nonzero coefficient, so $\AT(M)\le D+1$. Suppose now that $R_D=\{w\}$ is a single point, and put $\ell_i=\min(d_i,D)$. The points of $R_D$ are precisely the integer vectors $x$ with $0\le x_i\le \ell_i$ and $x_1+x_2+x_3=e$. Since $d_i>0$ and $D\ge D^*\ge1$, each $\ell_i$ is positive. If some $w_i<\ell_i$, choose $j\ne i$. If $w_j>0$, transfer one unit from coordinate $j$ to coordinate $i$; if $w_j=0$ and $w_i>0$, transfer one unit from coordinate $i$ to coordinate $j$; and if $w_i=w_j=0$, use the positive third coordinate as the donor and coordinate $i$ as the recipient. In each case this gives a second point of $R_D$, a contradiction. Hence $w_i=\ell_i=\min(d_i,D)$ for every $i$.

If $w_k=d_k$ for some $k$, say $k=3$, the exponent $d_3=\mu_{13}+\mu_{23}$ forces us to take the term $-u_3$ from every factor involving $u_3$. The remaining factor is $(u_1-u_2)^{\mu_{12}}$, and $w_1+w_2=e-d_3=\mu_{12}$. Thus $p_w=\pm\binom{\mu_{12}}{w_1}\ne0$. Otherwise $w=(D,D,D)$ with $3D=e$ and $D<\min_kd_k$. This forces $D=D^*$: if $D>D^*$, a point of $R_{D^*}$ would give a second point of $R_D$. If $p_w=0$, then $R_{D^*+1}$ contains the adjacent points $(D^*+1,D^*,D^*-1)$ and $(D^*+1,D^*-1,D^*)$. Proposition~\ref{prop:noadj}(b) guarantees that one of their coefficients is nonzero, and therefore $\AT(M)=D^*+2$.

For the second statement we compare \eqref{eq:twindensity} with the conditions defining $D^*$. For a twin pair $\{i,j\}$ the condition $2D-1\ge\mu_{ij}$ is equivalent to $2D\ge\mu_{ij}$, because $\mu_{ij}$ is odd. For $J=\{1,2,3\}$ the condition is $3D\ge e$ if $M$ has no twins, $3D-1\ge e$ if it has exactly one twin pair, and $3D-3\ge e$ if all three vertices are mutual twins. In the last case all three multiplicities are equal, so $3D^*=e$. Hence the twin corrections matter only when $3D^*=e$ and $M$ has twins; if $3D^*>e$, the condition $3D-1\ge e$ already holds at $D=D^*$.

Suppose that $3D^*=e$ and $M$ has twins. Then \eqref{eq:twindensity} fails at $D^*$ and holds at $D^*+1$, so $D_0=D^*+1$. The central vector is fixed by a twin swap, so $p_{(D^*,D^*,D^*)}=0$ by the argument of Lemma~\ref{lem:twin}. Since $3D^*=e$, the set $R_{D^*}$ consists of the central vector alone, and it is not a vector with some $w_k=d_k$, because such vectors have nonzero coefficient. So $D^*<\min_kd_k$, and the first statement gives $\AT(M)=D^*+2=D_0+1$. In all other cases \eqref{eq:twindensity} holds at $D^*$. Moreover $D^*\ge2=\chi-1$, since $D^*=1$ forces $e=3$, that is $\mu=(1,1,1)$, which has twins and $3D^*=e$. Hence $D_0=D^*$, and the second statement follows from the first.
\end{proof}

\begin{remark}
Recurrences of this kind are classical in another guise: up to normalisation, the coefficients of $P$ are the Wigner $3j$ symbols of quantum mechanics, and Lemma~\ref{lem:rec} is equivalent to a recurrence of Schulten and Gordon \cite{SG}. Zeros of these symbols that are not explained by symmetry, such as the one in Example~\ref{ex:3jzero} below, have also been studied there \cite{Pain}. None of this is used in the proofs.
\end{remark}

\begin{corollary}\label{cor:three}
\leavevmode
\begin{enumerate}[label=\textup{(\alph*)}]
\item Every uniquely $3$-colourable graph $G$ satisfies
\[
D_0(M_G)-2\le f(G)\le D_0(M_G)-1.
\]
\item Conjecture~\ref{conj:M} holds for $\vec m=(m_1,m_2,m_3)$ unless $M^+(\vec m)$ is \emph{central-tight}, that is, it has no twins, $3\mid e$, $D_0=e/3$ and $D_0<\min_kd_k$. In the central-tight case the conjecture for $\vec m$ is equivalent to
\[
\sum_s(-1)^s\binom{\mu_{12}}{s}\binom{\mu_{13}}{D_0-s}\binom{\mu_{23}}{D_0-\mu_{12}+s}\ne0 .
\]
\item For all $m\ge2$, $f(K_{2,2,m})=\lceil(2m-4)/3\rceil$.
\end{enumerate}
\end{corollary}

\begin{proof}
(a) and (b) follow from Theorems~\ref{thm:collapse} and \ref{thm:three}. For (c), $\mu=(2,m,m)$, there are no twins and $D_0=\lceil(2m+2)/3\rceil$. If $m\not\equiv2\pmod3$ then $3\nmid e$. If $m\equiv2\pmod3$ and $D=(2m+2)/3$, the central sum equals $2\bigl[\binom mD\binom m{D-2}-\binom m{D-1}^2\bigr]<0$, because $\binom mj^2>\binom m{j-1}\binom m{j+1}$ for $1\le j\le m$ (the ratio $\binom mj/\binom m{j-1}=(m-j+1)/j$ is strictly decreasing in $j$). In both cases $f=D_0+1-3$.
\end{proof}

For complete tripartite graphs, Corollary~\ref{cor:three}(a) leaves only two possible values of $f$. We have checked by computer that for every central-tight vector with $m_1,m_2,m_3\le60$ the central sum in Corollary~\ref{cor:three}(b) is nonzero modulo the prime $2^{61}-1$ or modulo $10^9+7$, hence nonzero, so by Corollary~\ref{cor:three}(b) Conjecture~\ref{conj:M} holds for all complete tripartite graphs with part sizes at most $60$. Whether this nonvanishing holds in general remains open (Question~\ref{q:tripartite}).

The vanishing in Theorem~\ref{thm:three} does occur for general multigraphs.

\begin{example}\label{ex:3jzero}
For $\mu=(2,9,13)$ we have $e=24$, $D^*=D_0=8$, degrees $(11,15,22)$, no twins, and the central coefficient is, up to sign,
\[
\sum_{s=0}^2(-1)^s\binom2s\binom9{8-s}\binom{13}{6+s}=15444-123552+108108=0 .
\]
Hence $\AT(M)=10=D_0+2$. This zero is not forced by any symmetry of $P$. Now let $H$ be a connected bipartite graph with parts of sizes $3$ and $4$ and $7$ edges, and let $G=H\vee\overline{K_5}$. Every proper $3$-colouring of $G$ uses two colours on $H$ and the third on the five added vertices, so $G$ is uniquely $3$-colourable, and $M_G$ has multiplicities $(2,9,13)$. By Theorem~\ref{thm:collapse}, $f(G)=7$, whereas the twin-density bound and the lower bound of Theorem~\ref{thm:main} both give $6$. (For the latter, a subgraph $G'$ with $a$ vertices in $H$ and $b$ added vertices has $e(G')-2|V(G')|\le e(H')+(a-2)(b-2)-4$, where $H'$ is the part of $G'$ in $H$; this is at most $7+15-4=18$ for $a=7$, and at most $6+12-4<18$ for $a\le6$, since $H$ is connected.) Thus Conjecture~\ref{conj:M} does not extend to all uniquely colourable graphs; the special form $\mu_{ij}=(m_i-1)(m_j-1)+1$ is essential.
\end{example}

For four or more vertices the analogue of Proposition~\ref{prop:noadj}(b) fails: a direct computation shows that for the four-vertex multigraph with $(\mu_{12},\mu_{13},\mu_{14},\mu_{23},\mu_{24},\mu_{34})=(2,5,1,5,5,4)$ the adjacent exponent vectors $(1,5,13,3)$ and $(1,5,14,2)$ both have coefficient $0$. Both lie in the Newton polytope of $P_M$ (the convex hull of the exponent vectors of its monomials), and this multigraph has no twins, so the vanishing is not forced by a symmetry.

\subsection{All parts but one of equal size}
The tool is Morris's constant term identity \cite{Morris}, which extends Dyson's identity with equal exponents. It is equivalent to Selberg's integral formula; see the survey of Forrester and Warnaar \cite{FW} for proofs and background. It states that for integers $a,b,\kappa\ge0$ and $\ell\ge1$,
\begin{equation}\label{eq:morris}
\CT\prod_{i=1}^\ell(1-x_i)^a\Bigl(1-\frac1{x_i}\Bigr)^b\prod_{1\le i\ne j\le\ell}\Bigl(1-\frac{x_i}{x_j}\Bigr)^\kappa=\prod_{j=0}^{\ell-1}\frac{(a+b+\kappa j)!\,(\kappa(j+1))!}{(a+\kappa j)!\,(b+\kappa j)!\,\kappa!}\ne0 .
\end{equation}

\begin{theorem}\label{thm:morris}
Let $\ell\ge1$, $p,q\ge1$ and $\vec m=(p,q,\dots,q)$ with $\ell$ parts of size $q$, so $\chi=\ell+1$. In $M^+(\vec m)$, the exceptional class is joined to every other class by $A=(p-1)(q-1)+1$ edges, and any two of the other classes are joined by $B=(q-1)^2+1$ edges. Put
\[
\kappa=\Bigl\lceil\frac B2\Bigr\rceil,\qquad T=A+(\ell-1)\kappa,\qquad D_1=\max\Bigl\{\Bigl\lceil\frac{\ell T}{\ell+1}\Bigr\rceil,\ (\ell-1)\kappa\Bigr\}.
\]
Then $\AT(M^+(\vec m))=D_0(\vec m)+1=D_1+1$. Consequently Conjecture~\ref{conj:M} holds for $\vec m$, and $f(K_{p,q,\dots,q})=D_1-\ell$.
\end{theorem}

The two terms in $D_1$ are the twin-density bounds of Proposition~\ref{prop:obstructions} for the set of all classes and for the set of the $\ell$ equal classes. For $p=q$ the theorem recovers Theorem~\ref{thm:equal}, and for $p=1$, $\ell=2$ it gives $f(K_{1,m,m})=\lfloor(m-1)^2/2\rfloor-1$ for $m\ge3$ (compare Example~\ref{ex:gaps}(b)). For $\ell=2$ it proves Conjecture~\ref{conj:M} for every complete tripartite graph with two equal parts; in particular, the central coefficient in Corollary~\ref{cor:three}(b) is nonzero whenever two of the three part sizes are equal.

\begin{proof}
Label the exceptional class by $0$ and the others by $1,\dots,\ell$. If $q$ is even then $B=2\kappa$; if $q$ is odd then $B=2\kappa-1$ is odd, and the classes $1,\dots,\ell$ are pairwise twins.

\emph{Lower bound.} Proposition~\ref{prop:obstructions} with $J=\{1,\dots,\ell\}$ gives $D\ge(\ell-1)\kappa$ in both parities, and with $J=\{0,\dots,\ell\}$ it gives $(\ell+1)D\ge\ell A+\ell(\ell-1)\kappa=\ell T$. Hence $D_1\le D_0\le\AT(M^+)-1$.

\emph{Upper bound.} We show that for every $0\le a\le A$ some nonzero coefficient of $P=P_{M^+}$ has largest exponent $\max\{\ell a,T-a\}$. Put $b=A-a$. Since $P$ is homogeneous, it suffices to work with
\[
P(1,x_1,\dots,x_\ell)=\prod_{i=1}^\ell(1-x_i)^A\prod_{1\le i<j\le\ell}(x_i-x_j)^B;
\]
the exponent of $x_0$ in a monomial is then determined by the degree. We also use $(1-x)^Ax^{-b}=\pm(1-x)^a(1-1/x)^b$.

If $q$ is even, let $\gamma=(\ell a,t,\dots,t)$ with $t=b+(\ell-1)\kappa$. Since $\Delta(x)^{2\kappa}=\pm\prod_{i\ne j}(1-x_i/x_j)^\kappa\prod_ix_i^{(\ell-1)\kappa}$ by Lemma~\ref{lem:schur}, the coefficient $[x^\gamma]P$ equals, up to sign, the constant term \eqref{eq:morris}.

If $q$ is odd, write $P=\Delta(x)Q$ with $Q=\prod_i(x_0-x_i)^A\Delta(x)^{2\kappa-2}$ symmetric in $x_1,\dots,x_\ell$. Let $t=b+(\ell-1)(\kappa-1)$ and let $\gamma$ have $x_0$-exponent $\ell a$ and remaining exponents $(t^\ell)+\delta_\ell$. By Lemma~\ref{lem:schur},
\[
[x^\gamma]P=[s_{(t^\ell)}]Q(1,x)=\frac1{\ell!}\CT\Bigl[Q(1,x)\prod_ix_i^{-t}\prod_{i\ne j}\Bigl(1-\frac{x_i}{x_j}\Bigr)\Bigr],
\]
which is $\pm\frac1{\ell!}$ times \eqref{eq:morris}. Its largest exponent other than $\ell a$ is $t+\ell-1=b+(\ell-1)\kappa$.

In both cases the coefficient is nonzero and the largest exponent of $\gamma$ is $\max\{\ell a,T-a\}$.

\emph{Choice of $a$.} Let $a^*=T/(\ell+1)$. If $a^*\ge A$, equivalently $\ell A\le(\ell-1)\kappa$, take $a=A$; then $\ell a\le T-a$, and the maximum is $T-A=(\ell-1)\kappa$. Otherwise take $a=\lfloor a^*\rfloor$; since $\ell\lfloor a^*\rfloor\le T-\lfloor a^*\rfloor$, the maximum is $T-\lfloor a^*\rfloor=\lceil \ell T/(\ell+1)\rceil$. In either case all exponents are at most $D_1$, so $\AT(M^+)\le D_1+1$. The formula for $f$ follows from Theorem~\ref{thm:collapse}.
\end{proof}

Write $D_1(p,q,\ell)$ for the quantity in Theorem~\ref{thm:morris}. Since $K_{\vec m}\subseteq K_{m_i,q,\dots,q}$ with $q=\max_{j\ne i}m_j$, Lemma~\ref{lem:gap} gives an upper bound for every complete multipartite graph: for any index $i$,
\[
f(K_{\vec m})\le D_1(m_i,q,\chi-1)-(\chi-1).
\]

\section{Open problems}\label{sec:open}

\begin{question}\label{q:tripartite}
Prove Conjecture~\ref{conj:M}. For $\chi=3$ the general bound is already within one of the truth by Corollary~\ref{cor:three}(a); the remaining issue is to show that the central coefficients of Corollary~\ref{cor:three}(b) never vanish for central-tight complete-tripartite vectors. For $\chi\ge4$, the first open case is two part sizes, each occurring at least twice; the relevant constant terms then have three different exponents and are not covered by the known evaluations of Dyson and Morris, or by Aomoto's extension of Selberg's integral \cite{Aomoto,FW}. Is $\AT(M^+(\vec m))-D_0(\vec m)$ bounded by a function of $\chi$? By Corollary~\ref{cor:three} it is at most $1$ when $\chi=3$.
\end{question}

\begin{question}
For graphs with several $\chi$-colourings, Proposition~\ref{prop:reduction} expresses each coefficient of $P_{G\vee K_k}$ as a signed sum over the colour partitions of $G$, and different colourings may cancel. For example, let $G$ be $K_{2,2,4}$ minus a perfect matching between its two parts of size $2$. Then $G$ has two colour partitions, and in each of them the two small classes span only two edges, so \eqref{eq:mu} gives multiplicities $(0,4,4)$. The multigraph $(0,4,4)$ has Alon--Tarsi number $4$, which would suggest $f(G)=1$; in fact the contributions of the two partitions cancel, and $f(G)=2$. For each colour partition $\pi$ of $G$ into $\chi$ classes, let $M_\pi$ be the multigraph defined by \eqref{eq:mu}, with multiplicities clipped below at $1$. In the example both clipped multigraphs are $(1,4,4)$, with Alon--Tarsi number $5$, and $\max_\pi(\AT(M_\pi)-\chi)=2=f(G)$. Is $|f(G)-\max_\pi(\AT(M_\pi)-\chi)|$ bounded by a function of $\chi$? (Example~\ref{ex:gaps} shows that the density bound of Proposition~\ref{prop:density} is not within an additive constant of $f$.)
\end{question}

\begin{question}\label{q:extremal}
By Corollary~\ref{cor:extremal}, $n^2/8-n/2\le F(n)\le n^2/8$. Is $F(n)=f(K_{\lceil n/2\rceil,\lfloor n/2\rfloor})$ for all $n$, that is, is the balanced complete bipartite graph exactly extremal?
\end{question}

\begin{question}
Is there a combinatorial proof of Theorem~\ref{thm:collapse}? Alon and Tarsi \cite{AT} showed that the coefficient of $\prod_vx_v^{d^+(v)}$ in $P_G$ equals, up to sign, the difference between the numbers of even and odd Eulerian subdigraphs (spanning subgraphs in which every vertex has equal in- and out-degree) of any orientation of $G$ with out-degrees $d^+(v)$. A combinatorial proof might therefore come from a sign-reversing involution on Eulerian subdigraphs of a suitable orientation of $G\vee K_k$. For equal parts, Zeilberger's proof of Dyson's identity \cite{Zeil}, which evaluates the constant term as a signed count of tournaments, suggests that such a proof might exist.
\end{question}

\begin{question}
Determine the order of magnitude of the paint analogue $\min\{k:\chi_P(K_{m,m}\vee K_k)=k+2\}$. By Theorem~\ref{thm:main} it is at most $\lfloor(m-1)^2/2\rfloor$. Is it linear in $m$?
\end{question}

\appendix
\section{A proof of Proposition~\ref{prop:reduction} via the Nullstellensatz}\label{app:qcn}

This appendix gives a second proof of Proposition~\ref{prop:reduction}, which applies Theorem~\ref{thm:qcn} directly to $P_{G\vee K_k}$ and makes the palette picture of Remark~\ref{rem:qcn} precise. Besides the tools of Section~\ref{sec:prelim}, it uses three further standard facts about symmetric functions \cite[Ch.~7]{EC2}.
\begin{enumerate}[label=\textup{(S\arabic*)},start=6]
\item Let $\omega\colon\Lambda\to\Lambda$ be the linear map with $\omega(s_\lambda)=s_{\lambda'}$. Then $\omega$ is a ring automorphism, $\omega(h_r)=e_r$, and $\omega(p_j)=(-1)^{j-1}p_j$.
\item The \emph{Hall inner product} is the bilinear form $\langle\cdot,\cdot\rangle$ on $\Lambda$ for which the $s_\lambda$ are orthonormal; thus $\langle\varphi,s_\lambda\rangle$ is the coefficient of $s_\lambda$ in $\varphi$. Since $\lambda\mapsto\lambda'$ is a bijection, $\omega$ preserves $\langle\cdot,\cdot\rangle$. By (S2), if $\ell(\lambda)\le r$, then $\langle\varphi,s_\lambda\rangle$ equals the coefficient $[s_\lambda]$ of the restriction of $\varphi$ to $r$ variables.
\item \emph{(Newton's identities.)} Every symmetric polynomial in $r$ variables is a polynomial in the power sums $p_1,\dots,p_r$.
\end{enumerate}

As in Section~\ref{sec:reduction}, fix $k$ and an admissible monomial $x^\alpha y^b$ with $b=\delta_k+\lambda$ and deficits $r_v=D-\alpha_v$. In addition, choose for each vertex $v$ an arbitrary set $Z_v\subseteq[0,D]$ with $|Z_v|=r_v$, and for a set $S$ of $\chi$ numbers put
\[
\Phi_{G,Z}(S)=\sum_{c\colon V(G)\to S}\frac{P_G(c)\prod_v\prod_{z\in Z_v}(c_v-z)}{\prod_vN_S(c_v)}.
\]
This is the colouring sum \eqref{eq:Phi} with $c_v^{\,r_v}$ replaced by the polynomial $\prod_{z\in Z_v}(c_v-z)$, which has the same leading term.

\begin{lemma}\label{lem:appZ}
There is a symmetric polynomial $\Phi_{G,Z}(u_1,\dots,u_\chi)$ of degree at most $|\lambda|$ whose value at every $\chi$-element set $S$ is $\Phi_{G,Z}(S)$. Its homogeneous component of degree $|\lambda|$ is the polynomial $\Phi_{G,r}$ of Lemma~\ref{lem:Phipoly}.
\end{lemma}

\begin{proof}
Write $F(x)=P_G(x)\prod_v\prod_{z\in Z_v}(x_v-z)=\sum_af_ax^a$. For each monomial the sum over $c\colon V(G)\to S$ splits into one sum per vertex, so by Lemma~\ref{lem:dd} with $|S|=\chi$,
\[
\Phi_{G,Z}(S)=\sum_af_a\prod_v\Bigl(\sum_{t\in S}\frac{t^{a_v}}{N_S(t)}\Bigr)=\sum_af_a\prod_vh_{a_v-\chi+1}(S).
\]
The right-hand side is a symmetric polynomial in the elements of $S$. Each term has degree $\sum_v(a_v-\chi+1)=|a|-n(\chi-1)$ or vanishes, and $|a|\le\deg F=e(G)+\sum_vr_v$. By the budget equation \eqref{eq:budget}, $\sum_vr_v=(\chi-1)n-e(G)+|\lambda|$, so every term has degree at most $|\lambda|$. Terms of degree exactly $|\lambda|$ come only from monomials of $F$ of top degree, and the top-degree part of $F$ is $P_G(x)\prod_vx_v^{r_v}$; the same computation applied to it gives $\Phi_{G,r}$, as in the proof of Lemma~\ref{lem:Phipoly}.
\end{proof}

It helps to think of $[0,D]$ as a palette of $D+1=\chi+k$ colours. In the grid sums below, the universal vertices take a set $Y$ of $k$ colours, and $G$ is coloured from the remaining set $S=[0,D]\setminus Y$ of exactly $\chi$ colours.

\begin{lemma}\label{lem:appgrid}
We have
\[
\coef{x^\alpha y^b}P_{G\vee K_k}=\sum_{y\in[0,b_1]\times\dots\times[0,b_k]}\frac{V(y)\,\Phi_{G,Z}\bigl([0,D]\setminus\{y_1,\dots,y_k\}\bigr)}{N_{[0,b_1]}(y_1)\cdots N_{[0,b_k]}(y_k)},
\]
where terms in which two $y_i$ coincide are zero. That is, the universal vertices choose the colours $y_1,\dots,y_k$, and $G$ is coloured from the remaining $\chi$ colours.
\end{lemma}

\begin{proof}
Apply Theorem~\ref{thm:qcn} to $P_{G\vee K_k}=P_G(x)\prod_v\prod_i(x_v-y_i)\,V(y)$, using the grid $[0,D]\setminus Z_v$ (which has $\alpha_v+1$ points) for $x_v$ and $[0,b_i]$ for $y_i$ (since $\lambda\subseteq(\chi^k)$, we have $b_i\le b_1=k-1+\lambda_1\le D$, so every $y_i$ lies in $[0,D]$); the degree condition holds with equality because the monomial is admissible. This gives
\[
\coef{x^\alpha y^b}P_{G\vee K_k}=\sum_{y}\frac{V(y)}{\prod_iN_{[0,b_i]}(y_i)}\;\sum_{x\in\prod_v([0,D]\setminus Z_v)}P_G(x)\prod_v\frac{\prod_i(x_v-y_i)}{N_{[0,D]\setminus Z_v}(x_v)} .
\]
If two $y_i$ coincide, then $V(y)=0$. So fix $y$ with distinct entries, let $Y=\{y_1,\dots,y_k\}$ and $S=[0,D]\setminus Y$, which has exactly $D+1-k=\chi$ elements, and evaluate the inner sum. Consider one vertex $v$ and a value $t=x_v\in[0,D]\setminus Z_v$. If $t\in Y$, then $\prod_i(t-y_i)=0$. If $t\in S\setminus Z_v$, then
\[
\frac{\prod_i(t-y_i)}{N_{[0,D]\setminus Z_v}(t)}
\overset{(1)}{=}\frac{\prod_{y\in Y}(t-y)\,\prod_{z\in Z_v}(t-z)}{\prod_{s\in[0,D],\,s\ne t}(t-s)}
\overset{(2)}{=}\frac{\prod_{y\in Y}(t-y)\,\prod_{z\in Z_v}(t-z)}{\prod_{y\in Y}(t-y)\cdot N_S(t)}
=\frac{\prod_{z\in Z_v}(t-z)}{N_S(t)} .
\]
Here (1) holds because $N_{[0,D]\setminus Z_v}(t)$ is the product of $t-s$ over all $s\ne t$ in $[0,D]$ except those in $Z_v$, all of which are different from $t$; and (2) holds because $[0,D]\setminus\{t\}=Y\sqcup(S\setminus\{t\})$. Thus the edges to the universal vertices cancel exactly the part of the weight that involves the universal colours. Consequently
\[
\sum_{x\in\prod_v([0,D]\setminus Z_v)}P_G(x)\prod_v\frac{\prod_i(x_v-y_i)}{N_{[0,D]\setminus Z_v}(x_v)}
=\sum_{x\in\prod_v(S\setminus Z_v)}\frac{P_G(x)\prod_v\prod_{z\in Z_v}(x_v-z)}{\prod_vN_S(x_v)}
=\Phi_{G,Z}(S).
\]
The first equality drops the vanishing terms with some $x_v\in Y$ and uses the computation above. The second holds because the summand vanishes whenever some $x_v\in Z_v$, so the sum may be extended to all $x\in S^{V(G)}$, which is the definition of $\Phi_{G,Z}(S)$.
\end{proof}

The sum in Lemma~\ref{lem:appgrid} still runs over the colours $y$ of the universal vertices, with $\Phi_{G,Z}$ evaluated at the complement of $\{y_1,\dots,y_k\}$. The next lemma evaluates such sums for an arbitrary symmetric polynomial. The idea is that a symmetric polynomial of the complement $[0,D]\setminus Y$ is a symmetric polynomial of $Y$, and $\omega$ translates between the two.

\begin{lemma}\label{lem:appomega}
Let $\lambda\subseteq(\chi^k)$ and $b=\delta_k+\lambda$. Let $\Phi$ be a symmetric polynomial in $\chi$ variables of degree at most $|\lambda|$, and let $\Phi_{|\lambda|}$ be its homogeneous component of degree $|\lambda|$. Then
\[
\sum_{y\in[0,b_1]\times\dots\times[0,b_k]}\frac{V(y)\,\Phi\bigl([0,D]\setminus\{y_1,\dots,y_k\}\bigr)}{N_{[0,b_1]}(y_1)\cdots N_{[0,b_k]}(y_k)}=(-1)^{|\lambda|}\,[s_{\lambda'}]\,\Phi_{|\lambda|}.
\]
Both sides are numbers: on the left $\Phi$ is evaluated at $\chi$-element sets, and on the right $[s_{\lambda'}]$ is a coefficient in the Schur expansion of $\Phi_{|\lambda|}$, as in Section~\ref{sec:prelim}.
\end{lemma}

\begin{proof}
\emph{Step 1: the left-hand side is a coefficient.} By Newton's identities, every symmetric polynomial in $\chi$ variables is a polynomial in the power sums $p_1,\dots,p_\chi$. So $\Phi(u)=\hat\Phi(p_1(u),\dots,p_\chi(u))$ for a polynomial $\hat\Phi$ in $\chi$ indeterminates. Give the $j$-th indeterminate weight $j$; we may take $\hat\Phi$ of weighted degree at most $|\lambda|$, and then $\Phi_{|\lambda|}(u)=\hat\Phi_{|\lambda|}(p_1(u),\dots,p_\chi(u))$, where $\hat\Phi_{|\lambda|}$ is the part of $\hat\Phi$ of weighted degree exactly $|\lambda|$. Define a polynomial in $y=(y_1,\dots,y_k)$ by
\[
\Psi(y)=\hat\Phi\bigl(p_1([0,D])-p_1(y),\dots,p_\chi([0,D])-p_\chi(y)\bigr).
\]
For distinct $y_1,\dots,y_k\in[0,D]$ we have $p_j([0,D]\setminus\{y_1,\dots,y_k\})=p_j([0,D])-p_j(y)$, so $\Psi(y)=\Phi([0,D]\setminus\{y_1,\dots,y_k\})$; terms with coinciding $y_i$ vanish on both sides because of $V(y)$. Hence the left-hand side equals $\sum_yV(y)\Psi(y)/\prod_iN_{[0,b_i]}(y_i)$. Since $\deg(V\Psi)\le\binom k2+|\lambda|=|b|$, Theorem~\ref{thm:qcn} identifies it with $\coef{y^b}(V\Psi)$. The monomial $y^b$ has the largest possible degree $|b|$, so only the part of $\Psi$ of degree exactly $|\lambda|$ contributes, namely
\[
\Psi_{\rm top}(y)=\hat\Phi_{|\lambda|}\bigl(-p_1(y),\dots,-p_\chi(y)\bigr),
\]
a symmetric polynomial in $y_1,\dots,y_k$.

\emph{Step 2: a Schur coefficient in $k$ variables.} By Lemma~\ref{lem:schur} with $r=k$, $\coef{y^{\lambda+\delta_k}}(V\Psi_{\rm top})=[s_\lambda]\Psi_{\rm top}$.

\emph{Step 3: passing to $\chi$ variables with $\omega$.} Let $F=\hat\Phi_{|\lambda|}(-p_1,\dots,-p_\chi)\in\Lambda$; its restriction to $k$ variables is $\Psi_{\rm top}$. Since $\lambda$ has at most $k$ parts, $[s_\lambda]\Psi_{\rm top}=\langle F,s_\lambda\rangle$, and since $\omega$ preserves the inner product, this equals $\langle\omega F,s_{\lambda'}\rangle$. By (S6), $\omega F=\hat\Phi_{|\lambda|}(\omega(-p_1),\dots,\omega(-p_\chi))$ with $\omega(-p_j)=(-1)^jp_j$. Restricted to $\chi$ variables $u$, $(-1)^jp_j(u)=p_j(-u)$, so the restriction of $\omega F$ is $\hat\Phi_{|\lambda|}(p_1(-u),\dots,p_\chi(-u))=\Phi_{|\lambda|}(-u)=(-1)^{|\lambda|}\Phi_{|\lambda|}(u)$. Finally $\lambda'$ has $\lambda_1\le\chi$ parts, so $\langle\omega F,s_{\lambda'}\rangle$ is the coefficient of $s_{\lambda'}$ in this restriction, which is $(-1)^{|\lambda|}[s_{\lambda'}]\Phi_{|\lambda|}$.
\end{proof}

\begin{proof}[Second proof of Proposition~\ref{prop:reduction}]
By Lemma~\ref{lem:appgrid}, the coefficient in question is the left-hand side of Lemma~\ref{lem:appomega} for $\Phi=\Phi_{G,Z}$. By Lemma~\ref{lem:appZ}, this $\Phi$ is a symmetric polynomial of degree at most $|\lambda|$ whose component of degree $|\lambda|$ is $\Phi_{G,r}$, so Lemma~\ref{lem:appomega} gives the claim.
\end{proof}


\end{document}